\documentclass{elsart}
\usepackage{amssymb}
\usepackage{graphicx}
\usepackage{latexsym}
\usepackage{amsfonts}
\usepackage[english]{babel}
\usepackage{xy}
\usepackage{fancyhdr}
\usepackage{latexsym,enumerate}
\usepackage{amsmath,amstext,amssymb}
\usepackage[dvips]{color}
\definecolor{rouge}{cmyk}{0,1,1,0.4}
\definecolor{bleut}{cmyk}{0.98,0.13,0,0.63}
\definecolor{bleur}{cmyk}{1,1,0,0.2}
\definecolor{bleuc}{cmyk}{0.98,0.13,0,0.43}
\definecolor{vert}{cmyk}{1,0,1,0.6}
\definecolor{gris}{gray}{.2}
\definecolor{grisc}{gray}{.3}
\usepackage{multicol}
\newcommand{\Z}{\ensuremath{\mathbb{Z}}}

\newcommand{\N}{\ensuremath{\mathbb{N}}}

\newtheorem{theorem}{Theorem}
\newtheorem{definition}{Definition}
\newtheorem{lemma}{Lemma}

\newtheorem{coro}{Corollary}[section]
\newtheorem{proposition}{Proposition}

\begin{document}

\begin{frontmatter}
\title{Thresholding methods to estimate
the copula density}
 \runtitle{Estimation of the copula density}

\author{F. Autin, E. Le Pennec, K. Tribouley}

\address{Universit\'es Aix-Marseille 1, Paris 7 et Paris 10}

\begin{abstract}\hspace{0,2cm}
This paper deals with the problem of the multivariate copula density
estimation. Using wavelet methods we provide two shrinkage
procedures based on thresholding rules for which the knowledge of
the regularity of the copula density to be estimated is not
necessary. These methods, said to be adaptive, are proved to perform
very well when adopting the minimax and the maxiset approaches.
Moreover we show that these procedures can be discriminated in the
maxiset sense. We produce an estimation algorithm whose qualities
are evaluated thanks some simulation. Last, we propose a real life
application for financial data.
\end{abstract}

\begin{keyword}
copula density, wavelet method, thresholding rules, minimax theory, maxiset theory.

\vspace{0.4cm}

\end{keyword}

\end{frontmatter}

\section{Introduction}
Recently a new tool has appeared (principally for risk management,
for instance in finance, insurance, climatology, hydrology etc) to
model the structure of dependance of the data. Let us start by
recalling the seminal result of Sklar \cite{sklar}.
\begin{theorem}[Sklar (1959)] Let $d\geq 2$ and $H$ be a $d-$variate
 distribution function.
 If, for $m=1,\ldots d$, any margin $F_m$  of $H$ is continuous, there exists a
  unique $d-$variate
 function $C$ with uniform margins ${\cal U}_{[0,1]}$ such that
\begin{displaymath}
\forall (x_1, \dots,x_d) \in \mathbb{R}^d,\;H(x_1,
\dots,x_d)=C(F_1(x_1),\dots, F_d(x_d)).
\end{displaymath}
\end{theorem}
The repartition function $C$ is called {\bf the copula } associated
to the distribution $H$. The interest of the Sklar Theorem is the
following: it is possible to study separately the laws of the
coordinates $X^m,m=1,\ldots d,$ of any vector $X$ whose the law is
$H$ and the dependance between the coordinates.

 The {\bf Copulas model} has been
extensively studied in  a parametrical framework for the
distribution function $C$. Large classes of Copulas, such as the
elliptic family, which contains the Gaussian Copulas and the Student
Copula, and the archimedian family, which contains the Gumbel
Copula, the Clayton Copula and the Frank Copulas, have been
identified. Mainly, people have worked in two directions. Firstly,
an important activity has concerned {\bf the modeling} in view to
find new Copulas and methodologies to simulate data coming from
these new Copulas. Secondly, usual {\bf statistical inference}
(estimation of the parameters, goodness-of-fit test, etc) has been
developed using the copulas. As usual, the nonparametric point of
view is useful when no a priori model of the phenomenon is
specified. For the practitioners, the nonparametrical estimators
could be seen as a {\bf benchmark} allowing to specify the model,
comparing with the available parametrical families. This explains
the success of the nonparametric estimator of the copula.
Unfortunately, the practitioners have difficulties to analyze the
graphes concerning the distribution functions (when $d=2$).
Generally, they try to make comments observing the scatter plot of
$\{(X_i,Y_i),i=1,\ldots,n\}$ (or $\{(R_i,S_i),i=1,\ldots,n\}$ where
$R,S$ are the rank statistics of $X,Y$). Therefore a density
estimator could be very advantageous.

We propose procedures to estimate the density associated to the
copula $C$. This density denoted $c$ is called the {\bf copula
density}. The copula density model that we consider is the
following.
 Let us give a $n$-sample $(X^1_1, \dots, X^d_1), \dots,$
$(X^1_n, \dots, X^d_n)$ of independent data admitting the same
distribution $H$ (and the same density $h$) as $(X^1,\dots,X^d)$. We
denote $F_1,\ldots,F_d$ the margins of the coordinates of the vector
$(X^1,\dots,X^d)$.  We are interested in the estimation of the
density copula $c$ which is defined as the derivative (if it exists)
of the copula distribution
\begin{displaymath}
c(u_1,\ldots,u_d)=\frac{h(F_1^{-1}(u_1), \dots,
F_d^{-1}(u_d))}{f_1(F_1^{-1}(u_1))\dots f_d(F_d^{-1}(u_d))}
\end{displaymath}
where
 $F_p^{-1}(u_p)=\inf\{x\in R:
F_p(x)\geq u_p\}$, $1 \leq p \leq d$ and  $u=(u_1, \dots,
u_d)\in[0,1]^d$.
 This model is a very classical density model but the direct
 observations
 $(U^1_i=F_1(X^1_i), \dots, U^d_i=F_d(X^d_i))$ for
$i=1,\ldots,n$ of the copula $C$ are not available, because the
margins are generally unknown. Observe that a similar problem could
be the nonparametric regression model with random design. This
model has been studied in Kerkyacharian and Picard \cite{KPwarped}
using {\bf warped wavelet families}.

In this paper, we focus on wavelet methods. We emphasize that these
methods are very well appropriated because the localization
properties of the wavelet allow us to give a sharp analysis of the
behavior of the copula density near the border of $[0,1]^d$. A great
advantage of the wavelet methods in statistics is to provide
adaptive procedures in the sense that they automatically adapt to
the regularity of the object to be estimated. The problem is a
little bit different in the case of the copula. So far, the
practitioners use copula densities which are regular. Nevertheless,
one of the practical problem is due to the fact that the copula
densities are generally concentrated towards the borders of
$[0,1]^d$ (see Part 6 where simulations are presented). It is very
difficult to provide a good estimation when the considered density
admits significant jumps (which is the case of the copula density
when its support is extended). Observe that the lack of data behind
the jumps produce poor estimators and this is just where we are
interested by the estimation. The thresholding methods use
multifrequency levels depending of the place where the estimator is
computed. Consequently, they are really fitted for the problem to
estimate the copula density and we think that they outperform the
classical linear methods like kernel estimators.

We present in this paper the {\bf Local Thresholding Method} and the
{\bf Global} (or {\bf Block}) {\bf Thresholding Method}. These
methods has been first studied by Donoho et al. (\cite{DJ2},
\cite{DJKP}) and Kerkyacharian et al. \cite{Kerk}. We prove that the
theoretical properties of these procedures with respect to the
quadratic risk are as good in the copula density model as in the
usual density model (when direct data are available). The good
behavior of the usual procedures in the copula density model is also
observed when linear procedures are considered (see Genest et al.
\cite{MasTrib}) and for test problems (see Gayraud and Tribouley
\cite{Gay}). To give an entire overview on the minimax properties of
the copula density estimators, we explore the maxiset properties of
both procedures (and also of the linear procedure). Again, we prove
that the fact that no direct observations are available does not
affect the properties of the wavelet procedures: as in the standard
density model, the local thresholding method of estimation
outperforms the others one.

 Next, we provide an important practical
study. Another advantage of the wavelet procedures is their
remarkable facility of use. Indeed, a large number of implementation
of the fast wavelet transform exist and it is just a matter of
choice to select one  in a given programming language. They can be
used almost directly after a simple warping of the original samples.
Nevertheless, as noticed previously in the case of the copula
density, as most of the pertinent information is located near the
border of $[0,1]^d$, the handling of the boundaries should be done
with a special care. It is one of the result of the article to show
that an inappropriate handling, such as the one proposed by default
by a lot of implementation, could cause the wavelet method to fail.
The symmetrization/periodization process proposed in this paper is
described. To reduce the gridding effect of the wavelet, we further
propose to incorporate some limited translation invariance in the
method. We first comment our results for simulated data in the case
of the usual parametrical copula families and then we present an
application on financial data: we propose a method to estimate the
copula density imposing that the "true" copula belongs to a target
parametrical family and using the nonparametrical estimator as a
benchmark.

 The paper is organized as follows. Section $2$ deals with the wavelet setting
so as to describe the multidimensional wavelet basis used in the
sequel. Section $3$ aims at describing both thresholding procedures
of estimation for which performances shall be studied in Section $4$
with the minimax approach and in Section $5$ with the maxiset
approach. Section $6$ is devoted to the practical results. Proofs of
main theorems are given in Section $7$ and shall need proposition
and technical lemmas proved in Appendix.
\section{Wavelet setting}
\noindent Let $\phi$ and $\psi$ be respectively a scaling function
and an associated wavelet function. We assume that these functions
are compactly supported on $[0,L]$ for some $L\geq1$. See for
example the Daubechies's wavelets (see Daubechies \cite{D}). For any
univariate function $h(\cdot)$, we denote by $h_{j,k}(\cdot)$ the
function $2^{j/2}h(2^j\cdot-k)$ where $j\in\N$ and $k\in \Z$.
 In the sequel, we use wavelet expansions for
multivariate functions. We build a multivariate wavelet basis as
follows: \begin{eqnarray*}
\phi_{j,k}(x_1,\dots,x_d)&=&\phi_{j,k_1}(x_1)\dots \phi_{j,k_d}(
x_d),\\
\psi_{j,k}^\epsilon(x_1,\dots,x_d)&=&\prod_{m=1}^d
\phi_{j,k_m}^{1-\epsilon_m}(x_m)\psi_{j,k_m}^{\epsilon_m}(x_m),
\end{eqnarray*}
for all $\epsilon=(\epsilon_1, \dots, \epsilon_d) \in S_d=\{0,1\}^d
\setminus \{(0,\dots,0)\}$. We keep the same notation for univariate
scaling function and multivariate scaling function; the subscripts
$k=(k_1,\dots, k_d)$ indicate the number of components. For any $j_0
\in \mathbb{N},$ the set $ \{\phi_{j_0k },\psi _{j,k}^\epsilon \;
|j\geq j_0,k\in
 \mathbb{Z}^d,\epsilon \in S_d\}$ is an orthonormal basis of
$L_2(\mathbb{R}^d)$ and the expansion of any real function $h$
of $L_2(\mathbb{R}^d)$ is given by:
\[\forall x\in \mathbb{R}^d,\;
h(x)=\sum_{k\in \mathbb{Z}^d}h_{j_0,k}\phi_{j_0,k}(x)+\sum_{j\geq
j_0}\sum_{k\in \mathbb{Z}^d}\sum_{ \epsilon \in S_d}h
_{j,k}^\epsilon \psi _{j,k}^\epsilon(x),
\]
where, for any $(j,k,\epsilon) \in \mathbb{N}\times
\mathbb{Z}^d\times S_d$, the wavelet coefficients are
\begin{eqnarray*} h_{j_0,k}=\int_{\mathbb{R}^d}h(x)
\phi_{j_0,k}(x)dx&\mbox{ and }& h_{j,k}^\epsilon =
\int_{\mathbb{R}^d}h(x) \psi_{j,k}^\epsilon(x)dx.
\end{eqnarray*}
Roughly speaking, the expansion of the analyzed function on the
wavelet basis splits into the
 "trend" at the level $j_0$ and the sum of the "details" for all the larger levels
$j,\;j\geq j_0$.

\section{Estimation procedures}
Assuming that the copula density $c$ belongs to $L_2$, we present
wavelet procedures of estimation. Motivated by the wavelet
expansion, we first estimate the coefficients of the copula density
on the wavelet basis. Observe that, for any $d-$variate function
$\Phi$
\begin{eqnarray*} E_c(\Phi(U_1,\ldots,U_{d}))&=
&E_h\left(\Phi(F_1(X^1),\ldots,F_{d}(X^{d}))\right)
\end{eqnarray*}
which means that the wavelet coefficients $
c_{j_0,k},c_{j,k}^\epsilon $ of the copula density $c$ on the
wavelet basis $ \{\phi _{j_0,k},\psi_{j,k}^\epsilon \; |j\geq
j_0,k\in \Z^{d},\epsilon \in S_d\}$ are equal to the coefficients of
the joint density $h$ on the warped wavelet family $$
\{\phi_{j_0,k}(F_1(\cdot),\ldots, F_{d}(\cdot)),\psi
_{j,k}^\epsilon(F_1(\cdot),\ldots, F_{d}(\cdot)) \; |j\geq j_0,k\in
\Z^{d},\epsilon \in S_d\}.$$ As usual, standard empirical
coefficients are
\begin{eqnarray}\label{coefemp}\widehat
{c_{j_0,k}}&=&\frac{1}{n}\sum_{i=1}^{n}\phi_{j_0,k}(F_1(X^1_i),\ldots,
F_{d}(X_i^{d}))\nonumber\end{eqnarray}
and
\begin{eqnarray}\label{coefemp}   \widehat
{c_{j,k}^\epsilon}&=&\frac{1}{n}\sum_{i=1}^{n}\psi_{j,k}^\epsilon(F_1(X^1_i),\ldots,
F_{d}(X_i^{d})).\end{eqnarray}
  Since no direct observation
$(F_1(X^1_i),\ldots, F_{d}(X_i^{d}))$ is usually available, we
propose to replace this one with the pseudo observation $(\widehat
{F_1}(X^1_i),\ldots, \widehat {F_{d}}(X_i^{d}))$ where $\widehat
{F_1},\ldots \widehat {F_{d}}$ are the empirical distribution
functions associated with the margins. The empirical coefficients
are
\begin{eqnarray*} \widetilde
{c_{j_0,k}}&=&\frac{1}{n}\sum_{i=1}^{n}\phi_{j_0,k}(\widehat
{F_1}(X^1_i),\ldots, \widehat {F_{d}}(X_i^{d}))=
\frac{1}{n}\sum_{i=1}^{n}\phi_{j_0,k}\left(\frac{R^1_i}{n},\ldots,\frac{R^{d}_i}{n}\right)
\end{eqnarray*}
and
\begin{eqnarray*} \widetilde {c_{j,k}^\epsilon}
&=&\frac{1}{n}\sum_{i=1}^{n}\psi_{j,k}^\epsilon(\widehat
{F_1}(X^1_i),\ldots, \widehat {F_{d}}(X_i^{d}))=
\frac{1}{n}\sum_{i=1}^{n}\psi_{j,k}^\epsilon\left(\frac{R^1_i}{n},\ldots,\frac{R^{d}_i}{n}\right)
\end{eqnarray*}
where $R^p_i$ denotes the rank of $X^p_i$ for $p=1,\ldots, d$
$$
R^p_i=\sum_{l=1}^n {\bf{1}}\{X^p_l<X^p_i\}.$$ According to the fact
that the wavelet basis is compactly supported, the sum over the
indices $k,\epsilon$ is finite and is taken over
$2^{d}L^{d}2^{{d}j}$ terms. In the sequel, to simplify the
notations, we omit the bounds of variation of the indices
$k,\epsilon$.

The most natural way to estimate the density $c$ is to reconstruct
the density thanks to the estimated coefficients. For any indices
$(j_n,J_n)$ such that $j_n \leq J_n$, we consider the very general
family of {\bf truncated estimators} of $c$ defined by
\begin{eqnarray}
\label{T} \widetilde {c_{T}}:=
\widetilde {c_{T}}(j_n,J_n) = \sum_{k}\widetilde
{c_{j_n,k}}\phi_{j_n,k}+\sum_{j=j_n}^{J_n}\sum_{k,\epsilon}\omega_{j,k}^\epsilon
\widetilde {c_{j,k}^\epsilon} \psi_{j,k}^\epsilon,
\end{eqnarray}
where $\omega_{j,k}^\epsilon \in \{0,1\}$ for any $(j,k,\epsilon)$.
It is intuitive that the more regular is the function $c$, the
smallest are the details and then that a good approximation for the
function $c$ is the trend at a level $j_n$ large enough. Such
{\bf{linear procedures}}
\begin{eqnarray}\label{L}
\widetilde
{c_{L}} := \widetilde {c_{L}}(j_n) = \sum_{k}\widetilde
{c_{j_n,k}}\phi_{j_n,k}
\end{eqnarray}
 consisting to put $\omega_{j,k}^\epsilon=0$  for any
$(j,k,\epsilon)$ have been considered in Genest et al.
\cite{MasTrib} where a discussion on the choice of $j_n$ is done.
 It is also possible to consider other
truncated procedures as {\bf non linear procedures}: the local
thresholding procedure consists to "kill" individually the small
estimated coefficients considering that they do not give any
information. Let $\lambda_n>0$ be a thresholding level and
$(j_n,J_n)$ be the level indices. For $\epsilon \in S_d$ and $j$
varying between $j_n$ and $J_n$, the hard local threshold estimators
of the wavelet coefficients are
\begin{eqnarray*}
\widetilde
{d_{j,k}^\epsilon}=\widetilde {c_{j,k}^\epsilon}
 {\bf{1}}\{|\widetilde {c^\epsilon_{j,k}}|> \lambda_n \},
 \end{eqnarray*}
  leading to the hard local threshold procedure
\begin{eqnarray}\label{HL}
\widetilde {c_{HL}}:=
\widetilde {c_{HL}}(j_n,J_n) =\sum_{k}\widetilde
{c_{j_n,k}}\phi_{j_n,k}+\sum_{j=j_n}^{
J_n}\sum_{k,\epsilon}\widetilde {d_{j,k}^\epsilon}
\psi_{j,k}^\epsilon.
\end{eqnarray}
It is also possible to decide to "kill" all the coefficients of the
level $j$ if the information given by this level is too small. For a
  thresholding level $\lambda_n>0$ and for levels $j$
varying between $j_n$ and $J_n$, we define the hard global threshold
estimates of the wavelet coefficients as follows
\begin{eqnarray*}
\widetilde
{e_{j,k}^\epsilon}=\widetilde {c_{j,k}^\epsilon} {\bf{1}}
\{\sum_k|\widetilde {c^\epsilon_{j,k}}|^2> L^{d} 2^{{d}j}\
\lambda_n^2 \},
\end{eqnarray*}
leading to the hard global threshold procedure
\begin{eqnarray}
\label{HG} \widetilde {c_{HG}}:= \widetilde
{c_{HG}}(j_n,J_n) =\sum_{k}\widetilde
{c_{j_n,k}}\phi_{j_n,k}+\sum_{j=j_n}^{
J_n}\sum_{k,\epsilon}\widetilde {e_{j,k}^\epsilon}
\psi_{j,k}^\epsilon.
\end{eqnarray}
The non linear procedures given in (\ref{HL}) and (\ref{HG}) depend
on the levels indices $(j_n,J_n)$ and on the thresholding level
$\lambda_n$ to be chosen by the user. In the next part, we explain
how to determine these parameters such a way that the associated
procedures achieve optimality properties. We need then to define a
criterion to measure the performance of our procedure.

\section{Minimax Results}
\subsection{Minimax approach}
A well-known way to analyze the performances of procedures of
estimation is the {\bf{minimax  theory}} which has been extensively
developed since the 1980-ies. In the minimax setting, the
practitioner chooses a loss function $\ell(.)$ to measure the loss
due to the studied procedure and a functional space $\mathcal{F}$
where the unknown object to estimated is supposed to belong. The
choice of this space is important because the first step (when the
minimax theory is applied) consists to compute the {\bf{minimax
risk}} associated to this functional space. This quantity
$$R_n(\mathcal{F})=\inf_{\widetilde {c}}\sup_{c \in
\mathcal{F}}E\,\ell(\widetilde {c}-c)$$(where the infimum is taken
over all the estimators of $c$) is a lower bound giving the best
rate achievable on the space $\mathcal{F}$. When a statistician
proposes an estimation procedure for functions belonging to
$\mathcal{F}$, he has to evaluate the risk of his procedure and
compare this upper bound with the minimax risk. If the rates
coincide, the procedure is {\bf{minimax optimal}} on the space
$\mathcal{F}$. A lot of minimax results for many statistical models
and many families of functional spaces as Sobolev spaces, Holder
spaces, and others as the family of Besov spaces have been now
established (see for instance
 Ibragimov and Khasminski \cite{IK} or Kerkyacharian and Picard \cite{KP}).

\subsection{Besov spaces}
Since we deal with wavelet methods, it is very natural to consider
Besov spaces as functional spaces because they are characterized in
term of wavelet coefficients as follows

\vspace{0.3cm}\begin{definition}[{\bf Strong Besov
spaces}]\label{Besov} For any $s>0$,
 a function $c$ belongs to the Besov space
$\mathcal{B}^s_{2\infty}$
 if and only if its sequence of wavelet coefficients $c_{j,k}^\epsilon$ satisfies
$$
\sup_{J\geq 0}2^{2Js}\sum_{j>
J}\sum_{k,\epsilon}(c_{j,k}^\epsilon)^2 < \infty.
$$
\end{definition}
 An advantage of the  Besov spaces $\mathcal{B}^s_{2
\infty}$ is to provide a useful tool to classify wavelet decomposed
signals according to their regularity and sparsity properties (see
for instance Donoho and Johnstone \cite{DJ2}). Last, it
is well known that the minimax risk measured with the quadratic loss
on this space is $$ \sup_n\inf_{\widetilde {c}}\sup_{c
\in \mathcal{B}^{s}_{2 \infty} }n^{\frac{2s}{2s+{d}}} \
E\,\|\widetilde {c}-c\|_2^2\ < \infty$$ where the infimum is taken
other any estimator of the density $c$.

\subsection{Optimality}
Let us focus on the quadratic loss function. Choosing a wavelet
regular enough and when ${d}=2$, Genest et al. \cite{MasTrib} prove
that the linear procedure $\widetilde {c_{L}}= \widetilde
{c_{L}}(j_n)$ defined in (\ref{L}) provides an
 optimal estimator on the Besov space $B^s_{2 \infty}$ for
 some fixed $s>0$ as soon as $j_n$ is chosen as follows:
$$ 2^{j_n -1}<  n^{\frac{1}{2s+{d}}}\leq 2^{j_n}.$$
On a practical point of view, this result is not completely
satisfying because the optimal procedure depends on the regularity
$s$ of the density which is generally unknown. To avoid this
drawback, many works in the nonparametric setting as in Cohen et al.
\cite{KP}, in Kerkyacharian and Picard \cite{KP3}
inspired from Donoho and Johnstone's studies on shrinkage procedures
(see for instance  \cite{DJ2}) build {\bf{adaptive
procedures}} which means that they do not depend on some a priori
information about the unknown density. For instance, note that the
thresholding procedures described in (\ref{HL}) and (\ref{HG}) are
clearly adaptive. The following theorem (which proof is a direct
consequence of Theorem \ref{maxih} established in the next section
by considering some inclusion spaces properties) gives results on
their rates
\vspace{0.3cm}\begin{theorem}\label{AdaptH} Let us consider a
continuously differentiable wavelet function and let $s>0$. Let us
choose the integers $j_n$ and $J_n$ and the real $\lambda_n$ such
that
$$2^{j_n-1}<\left(\log(n)\right)^{1/d}\leq 2^{j_n}, \quad 2^{J_n -1}<\left(\frac{n}{\log
n}\right)^{1/{d}}\leq 2^{J_n},
 \quad \lambda_n=\sqrt{\frac{\kappa\log(n)}{n}}$$
for some $\kappa>0$.
 Let $\widetilde {c}$ be either the hard local thresholding procedure $\widetilde {c_{HL}}(j_n,J_n)$
 or the hard global
  thresholding procedure
$\widetilde {c_{HG}}(j_n,J_n)$. Then, as soon as $\kappa$ is large
enough, we get
$$c \in \mathcal{B}^s_{2\infty} \cap L_\infty  \Longrightarrow  \sup_{n}
\left(\frac{n}{\log(n)}\right)^{\frac{2s}{2s+{d}}} E\| \widetilde
{c}-c\|_2^2< \infty. $$
\end{theorem}
We immediately deduce
\vspace{0.3cm}\begin{coro}\label{opti} The hard local thresholding
procedure $\widetilde {c_{HL}}$ and the hard global
  thresholding procedure
$\widetilde {c_{HG}}$ are adaptive near optimal up to a logarithmic
term (that is the price to pay for adaptation) on the Besov spaces
$\mathcal{B}^s_{2 \infty}$ for the quadratic loss function,
according to the minimax point of view.
\end{coro}
Observe that if $s>d/2$ then the assumption $c \in
\mathcal{B}^s_{2\infty} \cap L_\infty $ in Theorem \ref{AdaptH}
could be
 replaced  by $c \in \mathcal{B}^s_{2\infty}$ because
in this particular case $\mathcal{B}^s_{2\infty} \subsetneq L_\infty $.

\subsection{Criticism on the minimax point of view}On a practical point of view, the first
drawback of the minimax theory is the necessity to fix a functional
space. Corollary \ref{opti} establishes that no procedure could be
better on the space $\mathcal{B}^s_{2 \infty}$ than our hard
thresholding procedures (up to the logarithm factor) but this space
is generally an abstraction for the practitioner. Moreover, this one
knows that his optimal procedure
 achieves the best rate but he does not know this rate. An answer to this drawback is
  given by Lepski by introducing the concept of the random normalized factor
  (see Hoffmann and Lepski \cite{HL}).

 Secondly, he has the choice between both
procedures since
 Theorem \ref{AdaptH} establishes that the hard thresholding
procedures (local and global) have the same performances when
dealing with the minimax point of view.
 Nevertheless, a natural question arises here: is it possible to compare
thresholding procedures for estimation of the copula density between
themselves? To answer to these remarks, we propose to explore in the
next part the {\bf{maxiset approach}}.

\section{Maxiset Results}\label{maxiresults}
\subsection{Maxiset approach}
The maxiset point of view has been developed by Cohen et al.
\cite{CDKP} and is inspired from recent studies in the approximation
theory field. This approach aims at providing a new way to analyze
the performances of the estimation procedures. Contrary to the
minimax setting, the maxiset approach consists in finding the
maximal space of functions (called the maxiset) for which a given
procedure attains a prescribed rate of convergence. According to
this, the maxiset setting is not so far from the minimax one.
Nevertheless it seems to be more optimistic than the minimax point
of view in the sense that the maxiset approach points out {\bf
{all}} the functions {\bf {well-estimated}} by a fixed procedure at
a given accuracy.

In the sequel, we say that a functional space, namely
$\mathcal{MS}(\widetilde {c},r_n)$, is the maxiset of a fixed
estimation procedure $\widetilde {c}$ associated with the rate of
convergence $r_n$ if and only if the following equivalence is
satisfied
\begin{eqnarray*}
\sup _n  \ r_n^{-1} E\| \widetilde
{c}-c\|_2^2 < \infty \iff c
 \in \mathcal{MS}(\widetilde {c},r_n).
 \end{eqnarray*}
 Notice that the considered loss function is
 again the quadratic one.
 As a first consequence of adopting the maxiset point of
view, we  observe that if $\widetilde {c}$ is an estimator of
$c$ achieving the minimax rate of convergence $v_n$ on a
functional space $\mathcal{V},$ then $\mathcal{V}$ is included in
the maxiset of $\widetilde {c}$ associated with the rate
$v_n$ but is not necessarily the same. Therefore, it could be
possible to distinguish between both optimal minimax procedures: for
the same target rate, the best procedure is the procedure admitting
the largest maxiset.

 Recently many papers based on the maxiset approach have
arisen when considering the white noise model. For instance it has
been proved in Autin et al. \cite{AU3} that the hard local
thresholding procedure is the best one, in the maxiset sense, among
a large family of shrinkage procedures, called the {\bf {elitist
rules}},
 composed of all the wavelet procedures only using in their construction empirical wavelet
  coefficients with absolute value larger than a prescribed quantity. This
optimality has been already pointed out by Autin \cite{AU2} in
density estimation who has proved that {\bf {weak Besov spaces}} are
the saturation spaces of elitist procedures.

\subsection{Weak Besov spaces}
To model the sparsity property of functions, a very convenient and
natural tool consists in introducing the following particular class
of Lorentz spaces, namely {\bf {weak Besov spaces}}, that are in
addition directly connected to the estimation procedures considered
in this paper. We give definitions of the weak Besov spaces
depending on the wavelet basis. However, as is established in Meyer
\cite{M} and Cohen et al. \cite{CDKP}, most of them have also
different definitions proving that this dependence in the basis is
not crucial at all.

\vspace{0.3cm}\begin{definition}[{\bf Local weak Besov
spaces}]\label{WLBesov} \noindent For any $0<r<2$, a function $c \in
L_2([0,1]^d)$ belongs to the local weak Besov space
$\mathcal{W}_{_L}(r)$ if and only if its sequence of wavelet
coefficients $c_{j,k}^\epsilon$ satisfies the following equivalent
properties: \vspace{0.3cm}\begin{center}
\vspace{0.3cm}\begin{itemize}
\item $\vspace{0.3cm}\begin{displaystyle}\sup_{0<\lambda\leq 1}\lambda^{r-2}\sum_{j\geq
0}\sum_{k,\epsilon}
\end{displaystyle} (c_{j,k}^\epsilon)^2{\bf{1}}\{|c^\epsilon_{j,k}| \leq
\lambda \} < \infty,$ \item
$\vspace{0.3cm}\begin{displaystyle}\sup_{0<\lambda\leq
1}\lambda^{r}\sum_{j\geq 0}\sum_{k,\epsilon}
\end{displaystyle} {\bf{1}}\{|c^\epsilon_{j,k}| >
\lambda \} < \infty.$
\end{itemize}
\end{center}
\end{definition}

\vspace{0.3cm}

\begin{definition}[Global weak Besov
spaces]\label{WGBesov} \noindent For any $0<r<2$, a function $c \in
L_2([0,1]^d)$
belongs to the global weak Besov space $\mathcal{W}_{_G}(r)$ if and
only if its sequence of wavelet coefficients $c_{j,k}^\epsilon$
satisfies the following equivalent properties:
\vspace{0.3cm}\begin{center} \vspace{0.3cm}\begin{itemize}
\item $\vspace{0.3cm}\begin{displaystyle}\sup_{0<\lambda \leq 1}\lambda^{r-2}\sum_{j\geq
0}\sum_{k,\epsilon}
\end{displaystyle} (c_{j,k}^\epsilon)^2{\bf{1}}\{\sum_k (c^\epsilon_{j,k})^2 \leq
2^{{d}j} \lambda^2 \} < \infty,$
\item
$\vspace{0.3cm}\begin{displaystyle}\sup_{0<\lambda \leq
1}\lambda^{r}\sum_{j\geq 0}2^{{d}j} \sum_{\epsilon}
\end{displaystyle} {\bf{1}}\{\sum_k (c^\epsilon_{j,k})^2 >
2^{{d}j} \lambda^2 \} < \infty.$
\end{itemize}
\end{center}
\end{definition}
 The equivalences between the properties used in the
definition of global
 weak Besov spaces and in the definition of local
 weak Besov spaces can be proved
using same technics as those proposed in Cohen et al. \cite{CDKP}.

We prove in Section \ref{preuvelienbesov} the following link
 between the global weak Besov
space and the local weak Besov space
\vspace{0.3cm}\begin{proposition}\label{linclus} For any $0<r<2$, we
get $\quad \mathcal{W}_{G}(r) \subsetneq \mathcal{W}_{L}(r)$.
\end{proposition}

Last, we give an upper bound for the estimation expected error (when
the thresholding procedures are used) in the standard density model
(i.e. when direct observations $(U_i^1,\ldots,U_i^d),i=1,\ldots n$
are available). Notice that this result is stronger than the result
given in the minimax part because the functional assumption is
weaker.

\vspace{0.3cm}

\begin{theorem}\label{maxisetstandard} Let $\bullet$
be either $L$ or $G$ and  $\widehat{c_{H\bullet}}$ be the estimator
$\widehat{c_{H\bullet}}$ built in the same way as $\widetilde
{c_{H\bullet}}$ but with the sequence of coefficients $\widehat
{c_{j,k}^\epsilon}$ defined in (\ref{coefemp}).
 Let $s>0$ and assume that
\begin{eqnarray*}
c\in \mathcal{W}_\bullet(\frac{2{d}}{2s+{d}})\cap
 \mathcal{B}^{\frac{ds}{2s+ d}}_{2, \infty} .
\end{eqnarray*}
Then, there exists some $K>0$ such that for any $n$
\begin{eqnarray*}
E\|\widehat {c_{H\bullet}}-c\|_2^2&\leq& K \:
\left(\frac{\log(n)}{n}\right)^{\frac{2s}{2s+{d}}}
\end{eqnarray*}
\end{theorem}
There is no difficulty to prove this result by using the same
technics as in Cohen et al. \cite{CDKP} where the proof is done for
$d=1$  for classical thresholding rules in general models.
Nevertheless, for the interested reader, a detailed proof can be
found in Autin et al. \cite{ALT}.

\subsection{Performances and comparison of our procedures}\label{perf}
In this section, we study the maxiset performances of the linear
procedure and of the thresholding procedures described in Section
$1$. We focus on the optimal minimax procedures which means that we
consider the following choices of parameters:
$$2^{j_n^*-1}<\left(\frac{n}{\log(n)}\right)^{\frac{1}{2s+d}}\leq 2^{j_n^*},
 \quad \lambda_n=\sqrt{\frac{\kappa\log(n)}{n}}\mbox{ for some } \kappa>0,$$
$$2^{j_n-1}<\left({\log(n)} \right)^{1/{d}}\leq 2^{j_n},\quad
2^{J_n-1 }<\left(\frac{n}{\log(n)} \right)^{1/d} \leq 2^{J_n}.$$ Let
us fix $s>0$. we choose to focus on the (near) minimax rate
$r_n=\left(n^{-1}\log(n)\right)^{2s/(2s+{d})}$ achieved on the space
$B^s_{2\infty}$. The following theorem exhibits the maxisets of the
three procedures associated with the same target rate $r_n$.
\vspace{0.3cm}

\begin{theorem} \label{maxih}Let $s>0,$ and  assume
that $c\in L_\infty $. For a large choice of $\kappa$, we have
{\footnotesize
\begin{eqnarray}\label{maxih0}
&\hspace{1cm}&\sup _n
\left(\frac{n}{\log(n)}\right)^{\frac{2s}{2s+{d}}}
 E\|\widetilde {c_{L}}-c\|_2^2 < \infty \iff c
 \in \mathcal{B}^{s}_{2 \infty},\\
 \label{maxih1}&&
\sup _n \left(\frac{n}{\log(n)}\right)^{\frac{2s}{2s+{d}}}
 E\|\widetilde {c_{HL}}-c\|_2^2 < \infty \iff c
 \in \mathcal{B}^{\frac{ds}{2s+ d}}_{2 \infty} \cap
 \mathcal{W}_{L}(\frac{2d}{2s+ d}),\\
 &&\label{maxih2}
\sup _n \left(\frac{n}{\log(n)}\right)^{\frac{2s}{2s+{d}}}
 E\|\widetilde {c_{HG}}-c\|_2^2 < \infty \iff c
 \in \mathcal{B}^{\frac{ds}{2s+ {d}}}_{2 \infty} \cap
 \mathcal{W}_{G}(\frac{2d}{2s+ d}).
 \end{eqnarray}}
\end{theorem}
It is important to notice that, according to Theorem \ref{maxisetstandard}, the fact that direct observations are not available does not affect the maxiset performances of our procedures.
\noindent  As remarked by Autin et al. \cite{ALT} we have clearly $$\mathcal{B}^s_{2 \infty} \quad \subsetneq
\quad \mathcal{B}^{\frac{{d}s}{2s+{d}}}_{2 \infty}
  \cap
 \mathcal{W}_{G}(\frac{2d}{2s+d}).$$ We deduce by using Proposition \ref{linclus} that
both thresholding estimators considered in Theorem \ref{AdaptH}
 achieve the minimax rate (up to the logarithmic term) on a larger
 functional space than $\mathcal{B}^s_{2 \infty}$ which is the
 required space in the minimax approach. In particular, Theorem \ref{AdaptH} is proved.
 We propose now to discriminate these
procedures by comparing their maxisets. Thanks to Theorem
\ref{maxih} and applying the inclusion property given in Proposition
\ref{linclus}, we prove in Section \ref{preuvelienbesov} the
following corollary
\vspace{0.3cm}
\begin{coro}\label{bestmaxi} Let $s>0$ and let us
consider the target rate
$$r_n=\left(\frac{\log(n)}{n}\right)^{\frac{2s}{2s+{d}}}.$$ We get
\begin{eqnarray*}
{\mathcal{MS}}(\widetilde
{c_{L}},r_n) \subsetneq{\mathcal{MS}}(\widetilde {c_{HG}},r_n)
\subsetneq {\mathcal{MS}}(\widetilde {c_{HL}},r_n).
\end{eqnarray*}
 Hence, in the maxiset point of view and when the quadratic loss is considered, the
thresholding rules outperform the linear procedure. Moreover, the
hard local thresholding estimator $\widetilde {c_{HL}}$
appears to be the best estimator among the considered procedures
since it strictly outperforms the hard global thresholding estimator
$\widetilde {c_{HG}}$.
\end{coro}

\section{Applied results}\label{simulation}
This part is not only an illustration of the theoretical part.
First, we explain the considered algorithm with several numerical
possibilities to overcome any drawback. Next, we test the qualities
of our methodology with some simulation and we define the best
choices among our propositions. Last, we apply the chosen procedure
for financial data.
\subsection{Algorithms}

The estimation algorithms are described here for $d=2$ for sake of
simplicity but their extension in any other dimension is
straightforward. We therefore assume that a sequence
$\{(X_i,Y_i)\}_{1\leq i \leq n}$ of $n$ samples is given.

All estimators proposed in this paper can be summarized in a 7 steps
algorithm:
\begin{enumerate}
\item Rank the $X_i,Y_i$ with
\[
R_i = \sum_{l=1}^n \mathbf{1}_{X_l < X_i}\mbox{ and }S_i =
\sum_{l=1}^n \mathbf{1}_{Y_l < Y_i}.
\]
\item Compute the maximal scale index $J_n=\lfloor
  \frac{1}{2}\log_2(\frac{n}{\log n}) \rfloor$.
\item Compute the empirical scaling coefficients at the maximal scale index $J_n$:
\[
\widetilde{c_{J_n,k_1,k_2}} = \frac{1}{n} \sum_{i=1}^n
\phi_{J_n,k_1,k_2}(\frac{R_i}{n}, \frac{S_i}{n})\quad\text{for
  $1\leq k_1 \leq 2^{J_n}$ and $1\leq k_2 \leq 2^{J_n}$}.
\]
\item Compute the empirical wavelet coefficients
  $\widetilde{c^\epsilon_{j,k_1,k_2}}$ from this scaling
  coefficients with the fast 2D wavelet transform algorithm.
\item Threshold these coefficients according to the
  global thresholding rule or the local thresholding rule to obtained
  the estimated wavelet coefficients $\widetilde{e^\epsilon_{j,k_1,k_2}}$ or $\widetilde{d^\epsilon_{j,k_1,k_2}}$.
\item Compute the estimated scaling coefficients
  $\widetilde{c_{J_n,k_1,k_2}}$
at scale index $J_n$
  by the fast 2D wavelet inverse transform algorithm.
\item Construct the estimated copula density $\widetilde{c}$ by the formula
\[
\widetilde{c} = \sum_{k_1,k_2} \widetilde{c_{J_n,k_1,k_2}} \phi_{J_n,k_1,k_2}\quad.
\]
\end{enumerate}

Unfortunately only the steps $(1)$, $(2)$ and $(5)$ are as
straightforward as they seem. In all the other steps, one has to
tackle with two issues: the handling of the boundaries and the fact
that the results are not a function but a finite matrix of values.

The later issue is probably the easiest to solve. It means that we
should fix a number of point $N$ larger than $2^{J_n}$ and
approximate the estimated copula density at step $(7)$ on the
induced grid $(i/N,j/N)$. We can not compute exactly the value on
the grid as the scaling function is not always known explicitly.
Nevertheless a very good approximation can be computed on this grid
and we neglect the effect of this approximation. From the numerical
point of view, this implies that the estimation error can be
computed only on this grid and thus that the norms appearing in the
numerical results (see Table  \ref{RE_2_n500}, Table
\ref{RE_2_n2000}, Table \ref{RE_1_2_infty_n2000} and Table
\ref{RE_1_2_infty_n500}) are empirical norms on this grid. In this
paper, we choose $N=4*2^{J_n}$. Note that step $(3)$ also require an
evaluation of the scaling function and thus is replaced by an
approximation.

The former issue, the boundaries effect, is the key issue here.
Indeed, for most copula densities, the interesting behavior arises
in the corners which are the most difficult part to handle
numerically. The classical construction of the wavelet yields a
basis over $\mathbb{R}^2$ while we only have samples on $[0,1]^d$.
\begin{itemize}
\item A first choice is to consider the function  of $[0,1]^d$ to be
estimated as a function of $\mathbb{R}^2$ which is null outside
$[0,1]^d$. This choice is called {\bf zero padding} as it impose the
value $0$ outside $[0,1]^d$.
\item A second choice is to suppose that we observe the
restriction on $[0,1]^d$ of a $1$-periodic function, this is
equivalent to work in the classical periodic wavelet setting. This
choice called {\bf periodization} is very efficient if the function
is really periodic.
\item We propose also to modify the periodization and
assume that we observe the restriction over  $[0,1]^d$ of a even
$2$-periodic function. As we introduce a symmetrization over the
existing borders, we call this method {\bf symmetrization}. It
avoids the introduction of discontinuities along the border. Note
that nevertheless this symmetrization introduces a discontinuities
in the derivatives on the boundaries. \end{itemize}
 Once this
extension has been performed on the sample, one can apply the
classical wavelet transform. The resulting estimated copula density
is the restriction to $[0,1]^d$ of the estimated function.

The wavelet thresholding methods in a basis suffer from a gridding
effect: we can often recognize a single wavelet in the estimated
signal. To overcome this effect, we propose to use the {\bf cycle
  spinning} trick proposed by Donoho and Johnstone. To add some
translation invariance to the estimation process, we estimate the
copula density with a collection of basis obtained by a simple
translation from a single one and to average the resulting estimate.
In our numerical experiments, we have performed this operation with
$25$ different translation parameters and observed that it has
always ameliorate our estimate.
\subsection{Simulation}
We focus on the usual parametrical families of copulas: the FGM, the
Gaussian, the Student, the Clayton, the Frank and the Gumbel
families. We give results for both values of $n$ (the number of
data): $n=500$ is very small for bidimensional problems and $n=2000$
is usual in nonparametric estimation.

We test both methods of estimation and three ways to solve the
boundaries problems. We simulated data with the same copula, the
first margin being exponential with parameter $4$ and the second
margin being standard gaussian. Obviously, because of our algorithm,
the results are exactly the same when we change the laws of the
margins.

To evaluate the quality of our results, we consider three empirical
loss functions deriving from the $L_1$ norm, $L_2$ norm and
$L_\infty$ norm:
$$
 E_q=\|\widetilde{c}-c_0\|_{N,q}\quad \mbox{ for }
q=1,2,\infty
$$ where $c_0$ is the "true" copula density and $N\times N$ is the number
of points of the grid (see the previous part). Table
\ref{RE_1_2_infty_n2000} and Table \ref{RE_1_2_infty_n500}
 summarize the estimation relative errors
$$
 RE^q=\frac{\|\widetilde{c}-c_0\|_{N,q}}{\|c_0\|_{N,q}}\quad \mbox{ for }
q=1,2,\infty.
$$
These relative errors are computed with $100$ repetitions of the
experiment. The associated standard deviation is also given (in
parentheses).

\begin{figure}
  \centering

\begin{tabular}{ccc}
\includegraphics[width=7cm]{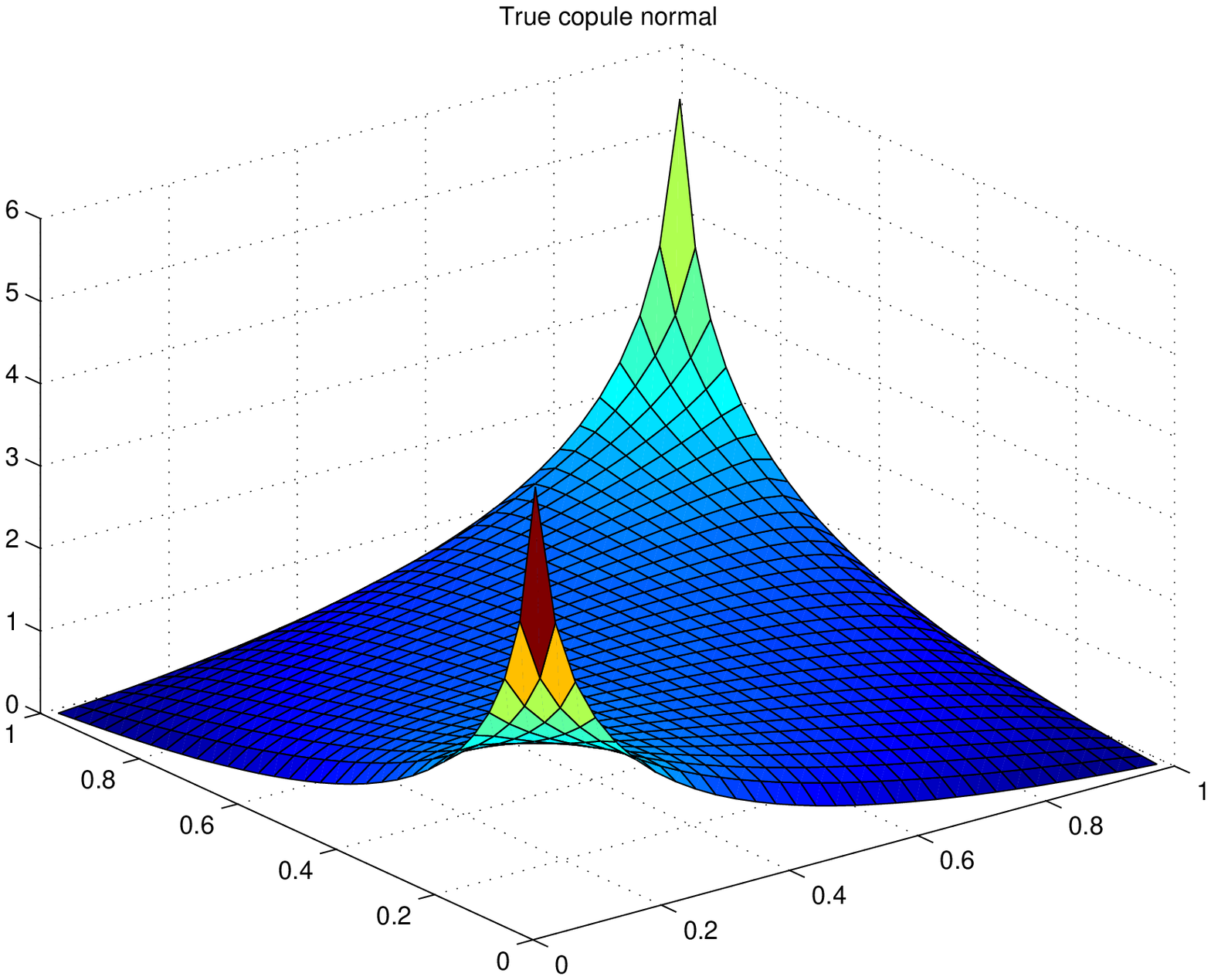}
&\hspace{1cm}
&\includegraphics[width=7cm]{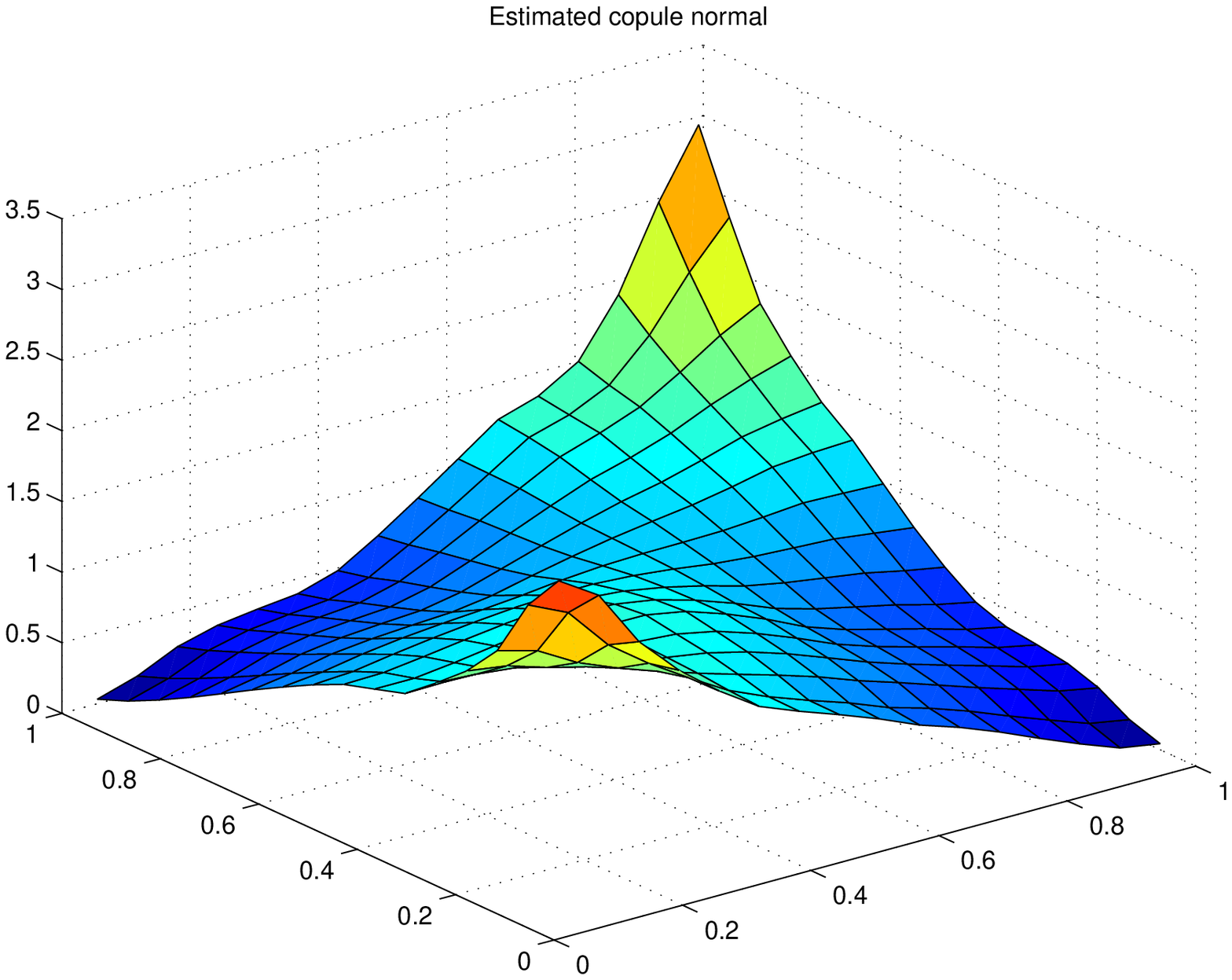}\\
\vspace{2cm}
(a) & &(b) \\
\includegraphics[width=7cm]{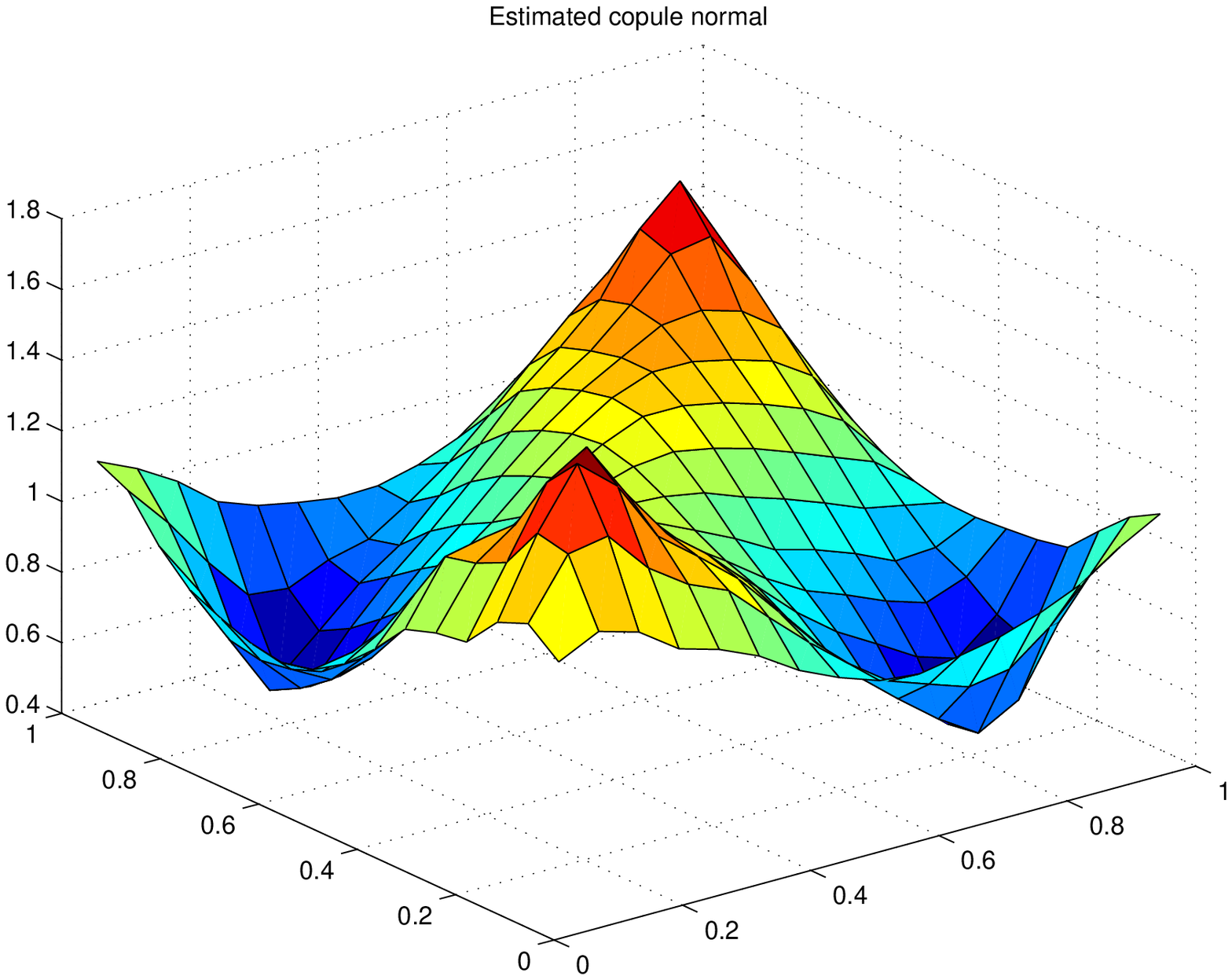}
& &
\includegraphics[width=7cm]{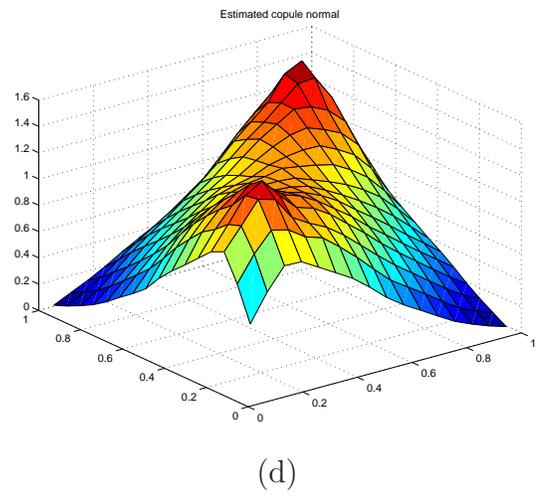}\\
(c) & &(d)
\end{tabular}
\caption{Estimation of the normal copula density of parameter $0.5$
with
  $n=2000$ (local thresholding): (a) true copula, (b) estimated copula with symmetrization, (c)
  estimated copula with periodization, (d) estimated copula with zero padding.}
\label{fig:extension}
\end{figure}

Table \ref{RE_2_n500} and  Table  \ref{RE_2_n2000} show that the
zero-padding method  and the periodization method give similar
results and lead to errors which are generally much larger than the
 symmetric periodization which is the best method to solve the
boundaries effects. This remark is valid for $n=500,2000$. Although
the zero-padding method is the default method with the Matlab
Wavelet Toolbox, it suffers from a severe drawback: it introduces a
strong discontinuities along the borders of $[0,1]^d$. The
periodization method suffers from the same drawback than the
zero-padding method as soon as the function is not really periodic.
Figure~\ref{fig:extension} emphasizes the superiority of the
symmetric periodization method in the case where the unknown copula
density is a normal copula. While the copula estimated with
symmetric extension remains close from the shape of the true copula
except for the resolution issue, this is not the case for the two
other estimated copulas: in the periodized version, the height of
the extreme peaks is reduced and two spurious peaks corresponding to
the periodization of the real peaks appear. The zero padded version
shows here only the reduced height artifact.

Tables~\ref{RE_1_2_infty_n500} and~\ref{RE_1_2_infty_n2000} display
the empirical $L_1$, $L_2$ and $L_\infty$ estimation error for the
symmetric extension for respectively $n=500$ and $n=2000$. Globally,
they show that the best results are obtained for the $L_2$ norm for
which the method has been designed. The second best results are
obtained for the $L_1$ norm because  a bound on the $L_2$ norm
implies a bound on the $L_1$ norm. The estimation in $L_\infty$ is
much harder as it is not a consequence of the estimation in $L_2$
and can be considered as a challenge for such a method.

One can also observed that the behavior largely depends on the copula
itself. This is coherent with the theory that states that the more
``regular'' the copula is the more efficient the estimator will
be. The copulas that are the least well estimated (Normal with
parameter .9, Student with parameter .5 and Gumbel with parameter
8.33) are the most ``irregular'' ones. They are very ``peaky'' for the
last two and almost singular along the first axis for the first one.
They are therefore not enough regular to be estimated correctly by the
proposed method.

A final remark should be given on the difficulty to evaluate such
these errors. Whereas the $L_1$ norm is finite equal to $1$ for all
true copula,
 the $L_2$ and $L_\infty$ norm can be very large even infinite because
 of their peaks. This is not an issue from the numerical point of view
 as we are restricted to a grid of step $1/N$ on which one can ensure
 the finiteness of the copula. Nevertheless the induce ``empirical''
 norm can be substantially different from the integrated norm. Thus
 the error  for $n=500$ to $n=2000$ are not strictly equivalent as the
 function can be much more complex  for the resolution induced by $n=2000$ than for $n=500$.

\subsection{Real data applications}
\begin{figure}[h]
\includegraphics[width=7cm,height=7cm]{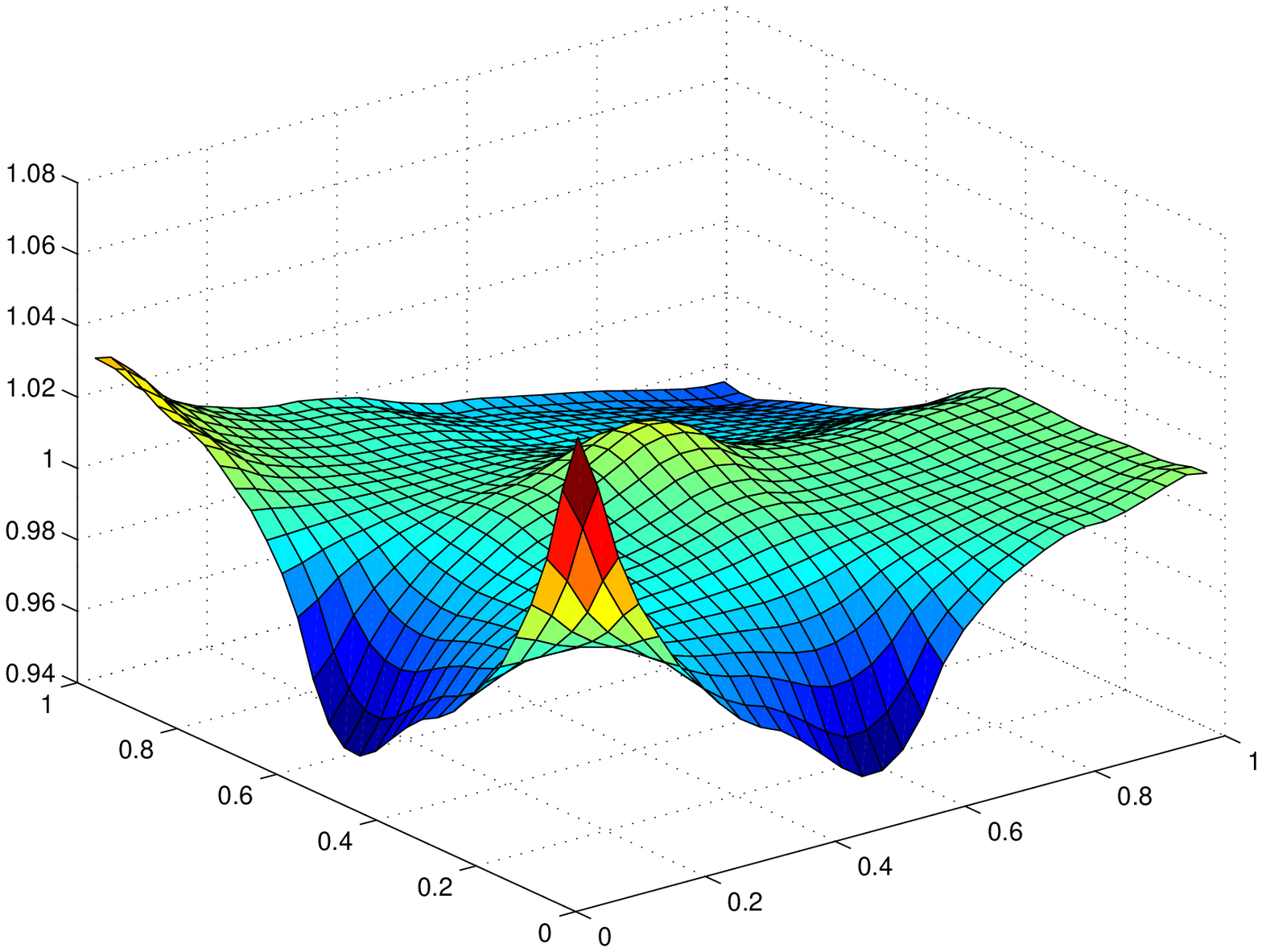}
\includegraphics[width=7cm,height=7cm]{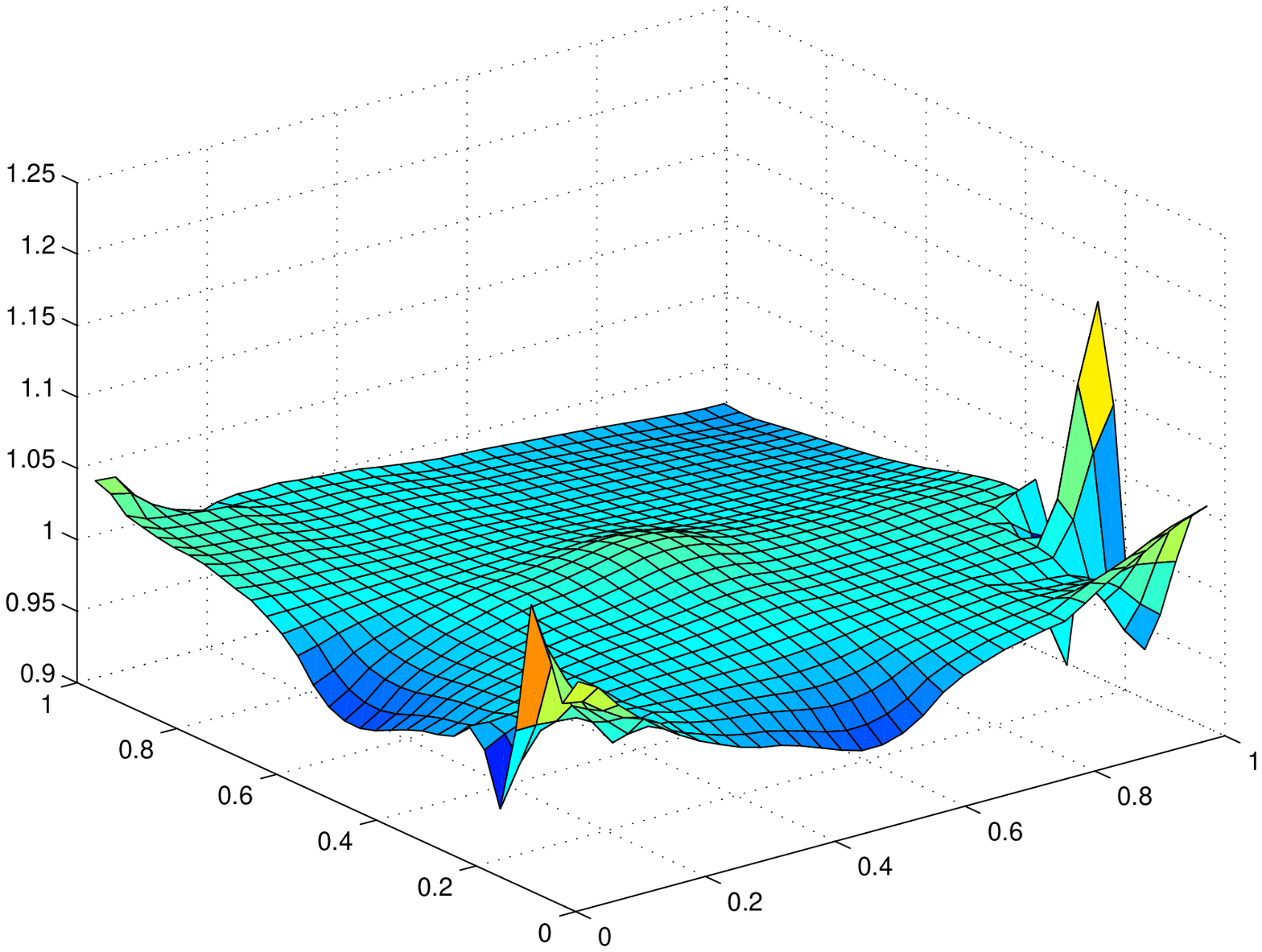}
\caption{{\small {Brent/Cac: Block Thresh. Method (left) and Local
Thresh. Method (right)}}\label{brent/cac}}
\end{figure}

\begin{figure}[h]
\includegraphics[width=7cm,height=7cm]{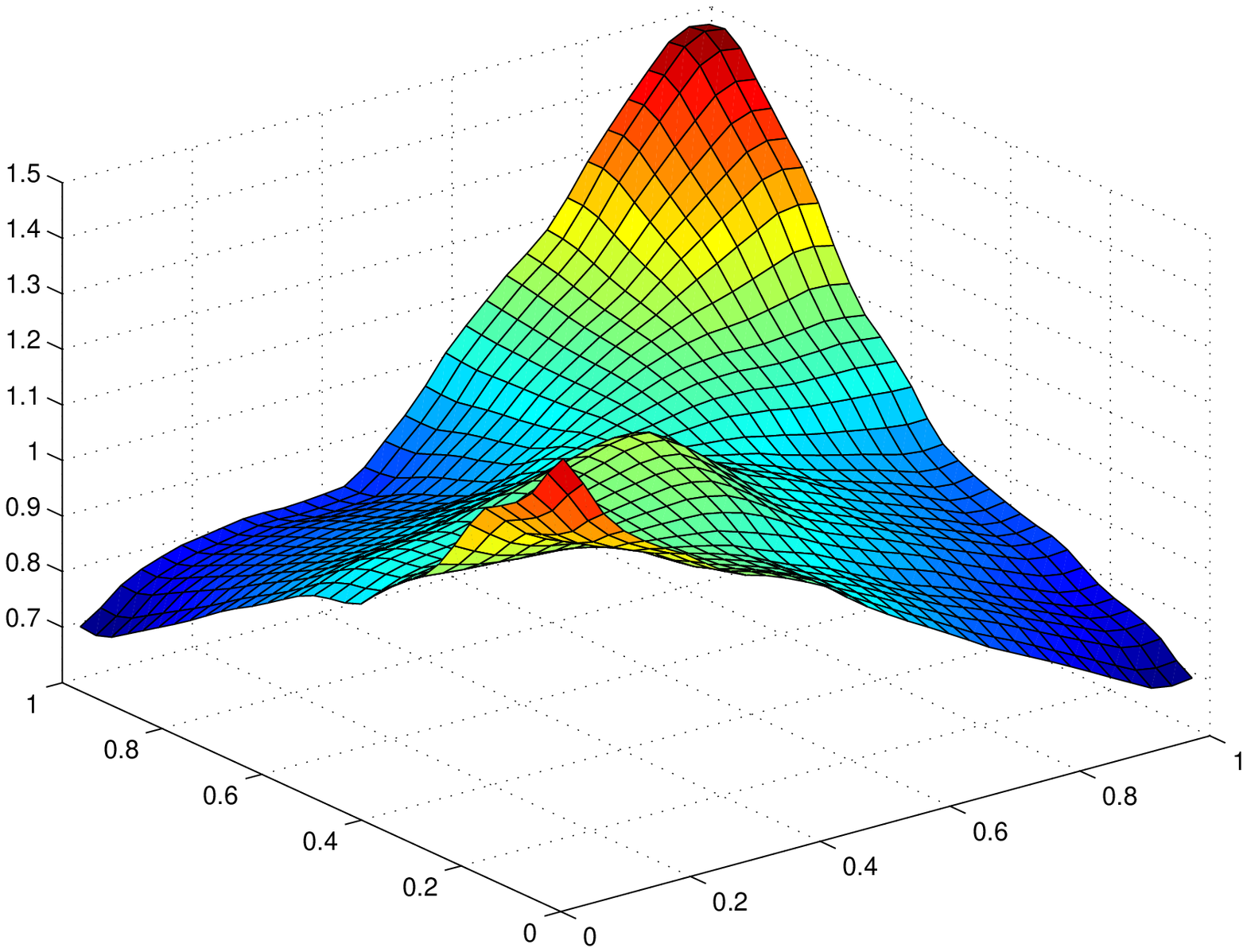}
\includegraphics[width=7cm,height=7cm]{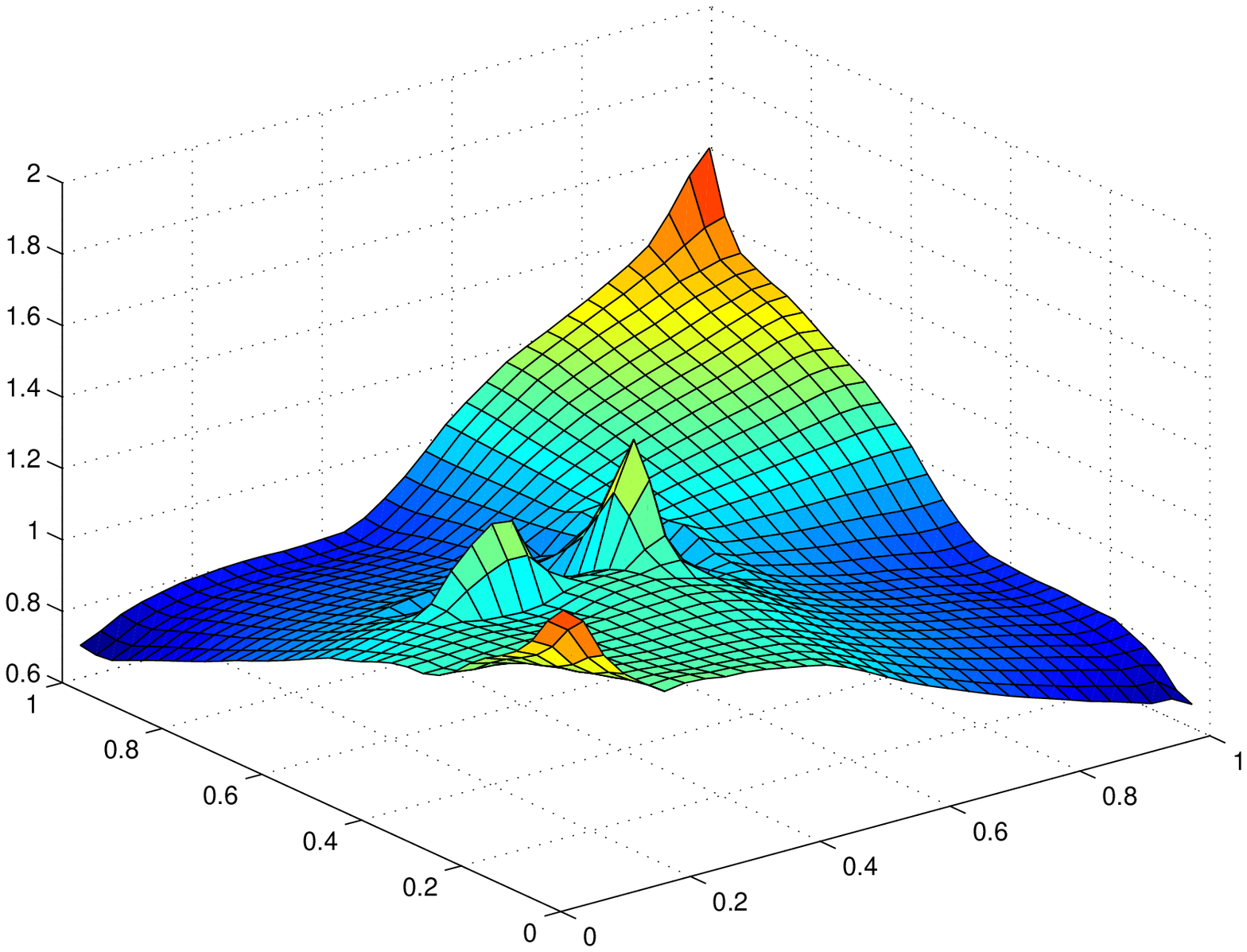}
\caption{{\small {Brent/ExonMobil: Block Thresh. Method (left) and
Local Thresh. Method (right)}}\label{brent/exonmobil}}
\end{figure}

\begin{figure}[h]
\includegraphics[width=7cm,height=7cm]{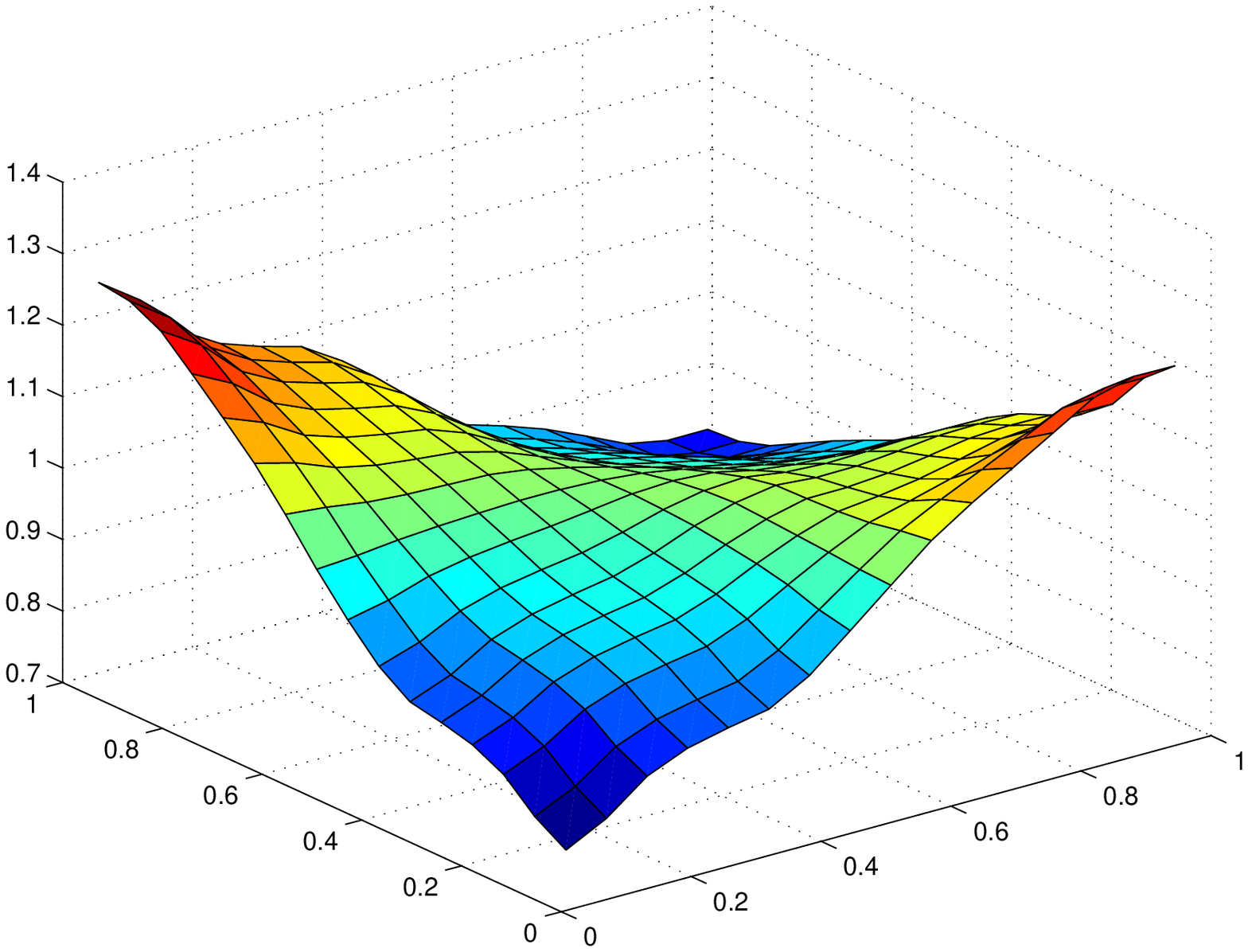}
\includegraphics[width=7cm,height=7cm]{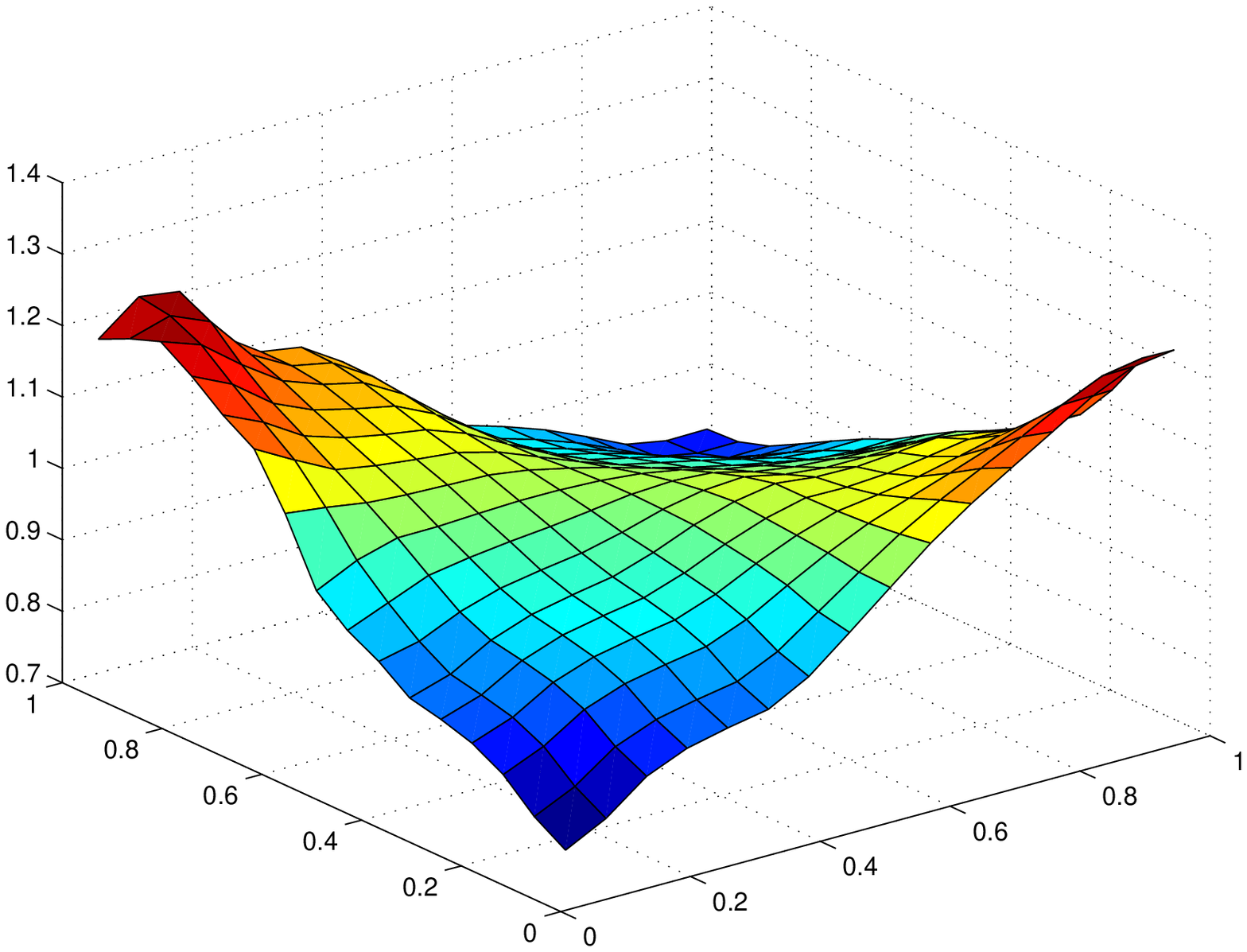}
\caption{{\small { DowJones/Oncedor: Block Thresh. Method (left) and
Local Thresh. Method (right)}}\label{dowjones/oncedor}}
\end{figure}

\begin{figure}[h]
\includegraphics[width=7cm,height=7cm]{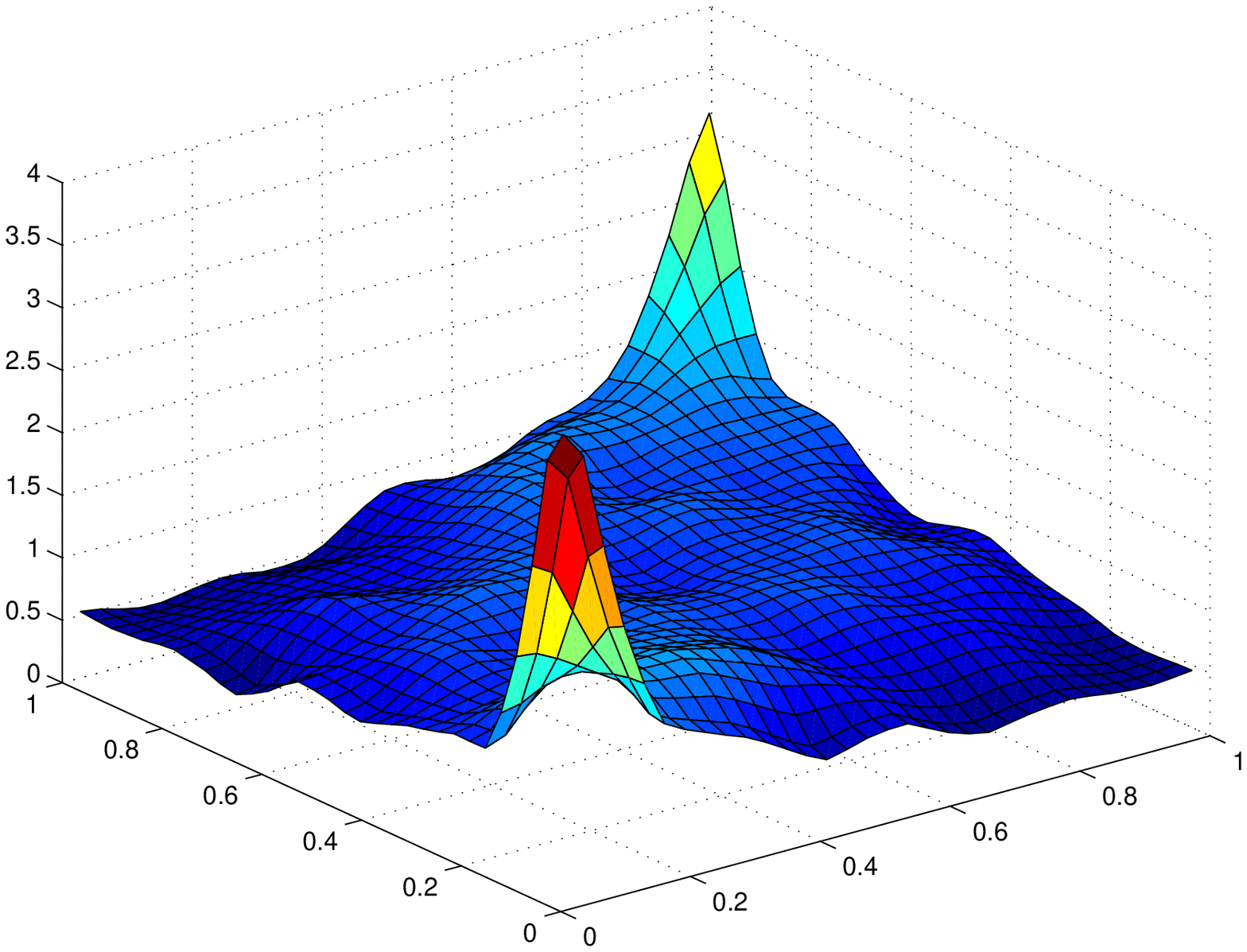}
\includegraphics[width=7cm,height=7cm]{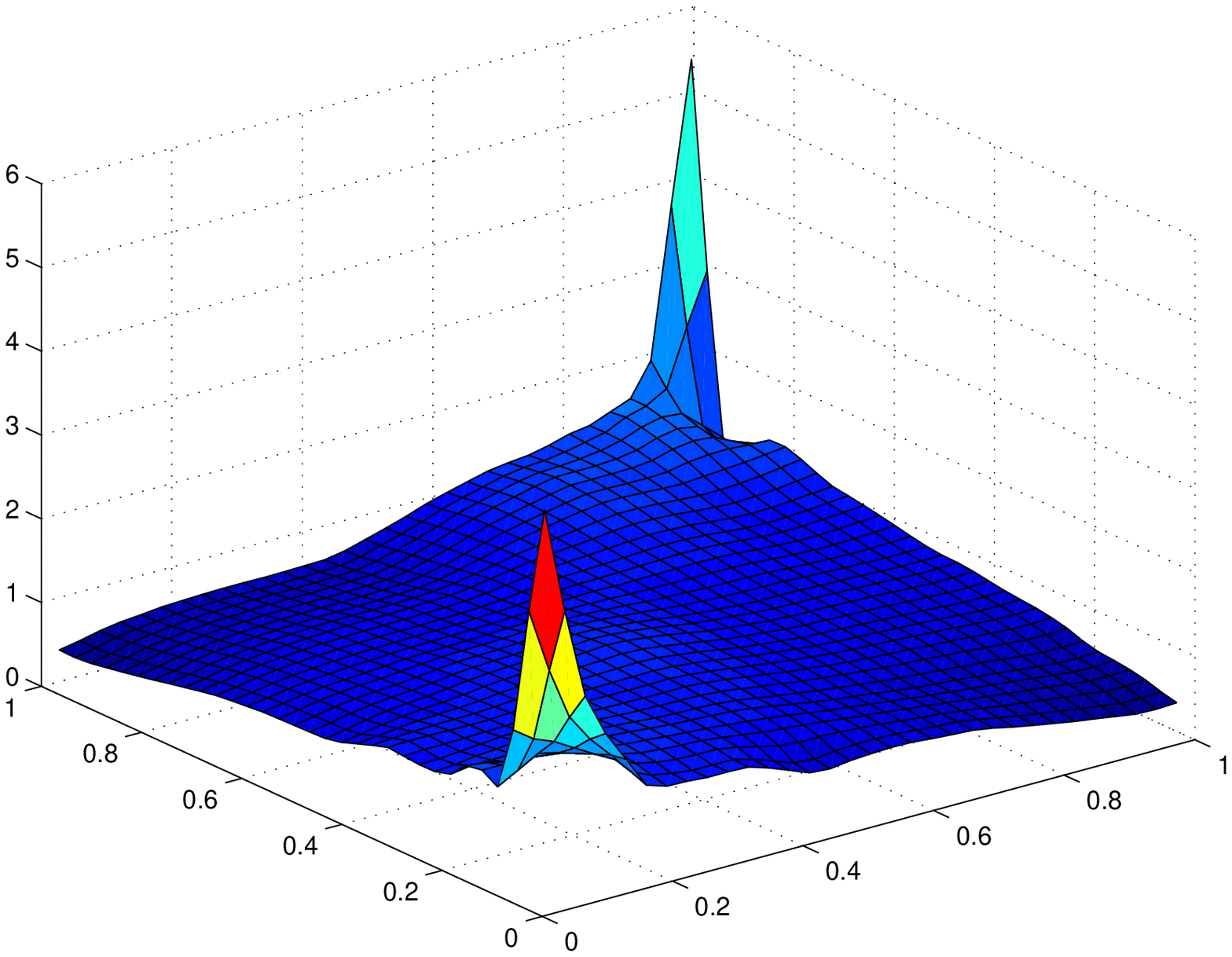}
\caption{{\small {DowJones/Fste100uk: Block Thresh. Method (left)
and Local Thresh. Method (right)}}\label{dowjones/ftse100uk}}
\end{figure}

We apply the thresholding methods on various financial series to
identify the behavior of the dependance (or non dependance). All our
data correspond to daily closing market quotations and are from
01/07/1987 to 31/01/2007. As usual, we consider the $\log-$return of
the data.
 Note that the data of each samples are not necessary independent but we apply our
procedures as there were. We first propose estimators of the
bivariate copula density associated with two financial series using
the adaptive thresholding procedures (see Figure \ref{brent/cac}--
Figure \ref{dowjones/ftse100uk}). Next, the nonparametrical
estimator denoted $\hat c$ is used as a benchmark and we derive an
new estimator by minimization  of the error between the benchmark
and a given parametrical family of copula densities.
 As previously, we focus on copulas which belong on the
Gaussian, Student, Gumbel, Clayton or Frank families.
 More precisely, we consider the following parametric
classes of copulas
\begin{eqnarray*}
{\mathcal C}_1&=&\left\{c\in {\mathcal N}_{\theta},\quad
\theta=[-0.99:0.01:0.99]\right\}\\
{\mathcal C}_2&=&\left\{c\in {\mathcal T}_{\theta} ,\quad
\theta=[-0.99:0.01:0.99, 1:1:100]\right\}\\
{\mathcal C}_3&=&\left\{c\in {\mathcal G}_{\theta},\quad
\theta=[1:0.01:2 ]\right\}\\
{\mathcal C}_4&=&\left\{c\in  {\mathcal C}_{\theta},\quad
\theta=[0:0.01:2 ]\right\}\\
{\mathcal C}_5&=&\left\{c\in  {\mathcal F}_{\theta},\quad
\theta=[-2:0.01:2 ]\right\}
\end{eqnarray*}
We consider three distances
\begin{eqnarray*}
E_q(\theta,p)&=&\|\hat c-c_\theta\|_{N,q} \quad \mbox{ for }\quad
q=1,2,\infty
\end{eqnarray*}
where $c_\theta\in {\cal C}_p,p=1,\ldots 5$. We propose to estimate
the parameter $\theta$ for each class ${\mathcal C}_p$ of copula
densities as follows
$$\hat \theta^q_p=\mathop{\mbox{arg min}}_{\theta}\,E_q(\theta,p)$$
which appears to be the best estimator of $\theta$ under the
constraint that the copula $c$ belongs to a fixed parametrical
family.
 We
derive estimators of $c$ among all the candidates $\left\{ c_{\hat
\theta^q_p},p=1,\ldots 5\right\}$ for each contrast $q=1,2,\infty$.
Table \ref{Table_brent/cac}-- Table \ref{Table_dowjones/ftse100uk}
give
\begin{itemize} \item the parameter $\tilde \theta^q$ for
$q=1,2,\infty$ defined by
$$
\tilde \theta^q=\mathop{\mbox{arg min}}_{p=1,\ldots,5}\left(
\mathop{\mbox{arg min}}_{\theta}\,E_q(\theta,p)\,\right),
$$
\item the parametric family ${\cal
C}_{\hat p}$ corresponding to the smallest error , \item the
associated relative errors defined by
$$
RE^q(\tilde\theta^q)=100\,\frac{\|\hat
c-c_{\tilde\theta^q}\|_{N,q}}{\|c_{\tilde\theta^q}\|_{N,q}}
$$
where $c$ is in ${\cal C}_{\hat p}$.
\end{itemize}
We test a lot of financial series and we select four revealing
examples: we never observe that the Clayton family or the Gumbel
family contains the best choice to model the unknown copula; the
families used are always the Student class or the Frank class.

First, we observe that the results are very good since any computed
relative error is small (in the worst case $RE^q\leq 20\%$). The
results are quite similar for both thresholding methods when the
unknown copula density does not present high picks (but the last
one: DowJones versus Ftse100uk). In a theoretical  point of view, we
prefer the block thresholding method because the estimator are
smoother. See by instance the case of the copula between Brent and
ExonMobil where the picks appearing in Figure \ref{brent/exonmobil}
(on the right) are not pleasant (even if their ranges are not so
large). Moreover, the relative error computed with the parametrical
density which is the best one among all the possible parametrical
copula densities is generally smallest for the block thresholding
method.

Notice that the choice of the contrast is crucial to estimate the
parameter $\theta$: there are significative differences between
$\hat \theta_1^p,\hat \theta_2^p,\hat \theta_\infty^p$. This is
usual in density estimation. We prefer to measure the loss due to
the estimator with the $L_1-$norm because this norm is very adapted
to the human eye and then the graphical results are the best. The
quadratic loss is frequently used because the graphical results are
easier to obtain but our opinion is that this norm does not
emphasize enough differences between the densities. See by instance
the very small relative errors computed with the $L_2$ contrast. The
$L_\infty$ norm has the opposite behavior: it accentuates every
difference. It could be a drawback when the local thresholding
method is considered and when to many details are kept (see again
the case of the copula of the couple Brent/ExonMobil).

Nevertheless, the choice of the best family do not depend on the
choice of the contrast:  it is fundamental because each type of
parametrical family is linked to a specific behavior of the
dependance and then the practitioner asks for indications about the
copula type. The study of the copula Cac versus Brent allows to
decide that both series are independent. Observe that there are
small problems on the borders (but notice the very small scaling)
although our methodology is made to remove this artefact. We think
that a usual linear kernel method becomes disastrous when the copula
density comes near of the uniform. The copula densities DowJones
versus Oncedor and Brent versus ExonMobil are both Frank copulas but
with opposite behaviors. It seems natural that the series Brent and
ExonMobil are dependent and varying in the same sense. Oncedor
(gold) is an hedge when the stock market collapses which could
explain the negative dependance between Oncedor and the financial
indices (we observe the same kind of dependance for others indices
like Fste100uk, Cac ...). The more delicate case is for the copula
DowJones versus Fste100uk because the picks are more accentuated. In
this case, the local thresholding method produces a nice Figure.

In conclusion, we present here an estimation method which is popular
among the practitioners: first, the nonparametrical estimator could
be see as a benchmark to decide graphically if the unknown copula
density looks like a copula density belonging to a well known
parametrical family. In this case, the parameter is estimated with
plug-in methods using the benchmark. We do not study here the
properties of such an estimator or the goodness-of-fit test problem.
For a statistic test procedure, we refer to Gayraud and Tribouley
(2008).

\section{Proofs}\label{proofs}
We first state the propositions needed to establish the main
results. Next, we prove Theorem \ref{maxih} in two steps by proving
both implications. Last, we establish Proposition \ref{linclus} and
Corollary \ref{bestmaxi}.

From now on we denote $K$ any constant that may change from one line
to another, which does not depend on $j$, $k$ and $n$ but which
depends on the wavelet and on $\|c\|_\infty$ and $\|c\|_2$.

\subsection{Preliminaries}\label{prel}
These preliminary results concern the estimation of the wavelet
coefficients and the scaling coefficients (denoted
$c_{j,k}^{\epsilon_0}$ with $\epsilon_0=(0,\ldots,0)$ to unify the
notation). Proposition \ref{variance} announces  that the accuracy
of the estimation is as sharp as if the direct observations were
available.

\vspace{0.3cm}

\begin{proposition}\label{LDL}
 Assume that the copula density belongs to $L_\infty$ and let
 $\delta>0$. There exists a constant $K>0$ such that for any $j$ such that
 $2^j\leq 2 \ \left(\frac{n}{\log(n)}\right)^{1/d}$,
and for any $(k, \epsilon)$
\vspace{0.3cm}\begin{eqnarray}
 \mathbb{P}\left(|\widetilde
{c^\epsilon_{j,k}}-\widehat {c^\epsilon_{j,k}}| > \lambda_n\right) &
\leq & K\,n^{-\delta}\label{a}\\
 \mathbb{P}\left(
 \sum_k(\widetilde{c^\epsilon_{j,k}}-\widehat{c^\epsilon_{j,k}})^{2} > L^{d}
\ 2^{{d}j} \lambda_n^2 \right) &\leq &
K\,n^{1-\delta}(\log(n))^{-1}\label{b}
\end{eqnarray}
as soon as $\kappa$ is chosen large enough.
\end{proposition}
It is clear that (\ref{b}) is a direct consequence of (\ref{a}).
Proof of (\ref{a}) is rejected to the Appendix. Note that from
(\ref{a}) we immediately deduce

\vspace{0.3cm}

\begin{proposition}\label{variance} Under the same
assumptions as in Proposition \ref{LDL} on $j$ and $c$, there exists a constant
$K>0$ such that for any $(k,\epsilon)$
\vspace{0.3cm}\begin{eqnarray*} E\left[\left(\widetilde
{c_{j,k}^\epsilon}-\widehat {c_{j,k}^\epsilon}\right)^2\right]
&\leq&K\,\frac{\log(n)}{n}.
\end{eqnarray*}
\end{proposition}

\subsection{Proof of Theorem \ref{maxih}}\label{ppl}
First, we prove the result for the linear estimator. Secondly, we
prove the result for the local thresholding method.  We do not prove
the result for the global thresholding method since the technics are
the same except that the required large deviation inequality
established in Proposition \ref{LDL} is given by (\ref{b}) instead
of (\ref{a}).

\subsubsection{Proof of Equivalence (\ref{maxih0})}
 On the one
hand, let $c$ be a copula density function belonging to $L_\infty $
satisfying for any $n$,
\begin{eqnarray}\label{hypo}
\mathbb{E}\|\widetilde {c_L}-c\|_2^2 \leq K
\left(\frac{\log(n)}{n}\right)^{\frac{2s}{2s+d}}
\end{eqnarray}
for some constant $K>0$. Let us prove that $c$ also belongs to the
space $ \mathcal{B}^s_{2\infty}$. Let us recall that the smoothing
index used for the linear procedure is $j_n^*$ satisfying
$2^{1-j_n^*} >  \left(n^{-1}\;\log(n)\right)^{1/(2s+d)}$. Since
\begin{eqnarray*}
 \mathbb{E}\|\widetilde
{c_L}-c\|_2^2&=&\mathbb{E}\|\widetilde {c_L}-\sum_k
c_{j_n^*,k}\phi_{j_n^*,k}\|_2^2+\|\sum_{j \geq
j_n^*}\sum_{k,\epsilon} c^\epsilon_{j,k}\psi^\epsilon_{j,k}\|_2^2,
\end{eqnarray*}
and following the assumption (\ref{hypo}),
 we get
 $$\sum_{j \geq
j_n^*}\sum_{k,\epsilon}(c_{j,k}^\epsilon)^2 \leq
\mathbb{E}\|\widetilde {c_L}-c\|_2^2 \leq K \ (2^{-2j_n^*})^{s}$$
which is the announced result. On the other hand, let us suppose
that $c\in \mathcal{B}^s_{2\infty}$. Then, using the
same technics as in Genest et al. \cite{MasTrib}, we prove that
$$\mathbb{E}\|\widetilde {c_L}-c\|_2^2 \leq K
\left(\frac{\log(n)}{n}\right)^{\frac{2s}{2s+d}}$$ which ends the
proof. The proof in Genest et al. \cite{MasTrib} is given in the
case $d=2$ and need some sharp control on the estimated coefficients
because an optimal result is established (there is no logarithmic
term in the rate).

\subsubsection{Proof  of Equivalence $(\ref{maxih1})$ (first step: $\Longrightarrow$)}

When direct observations $(F_1(X^1_i),\ldots,F_{d}(X^{d}_i))$ are
available, we use the estimator $\widehat  {c_{HL}}$ built in the
same way as $\widetilde {c_{HL}}$ but with the sequence of
coefficients $\widehat {c_{j,k}^\epsilon}$ defined in
(\ref{coefemp}) and the threshold $\lambda_n/2$ instead of
$\lambda_n$. Let us take $j_n,J_n$ positive integers and
$\lambda_n>0$. Since we get
\begin{eqnarray*} E\|\widetilde
{c_{HL}}-c\|_2^2&\leq & 2 E\|\widetilde {c_{HL}}-\widehat
{c_{HL}}\|_2^2+ 2 E\|\widehat  {c_{HL}}-c\|_2^2
\end{eqnarray*}
we have then to study the error term  due to the fact that we use
pseudo observations
 {\small \begin{eqnarray*}
T&=&E\|\widetilde {c_{HL}}-\widehat {c_{HL}}\|_2^2\\
&=&E \left[\sum_{k} (\widetilde {c_{j_nk}^{\epsilon_0}}-\widehat
{c_{j_nk}^{\epsilon_0}})^2\right] + E
\left[\sum_{j_n}^{J_n}\sum_{k,\epsilon} (\widetilde
{c_{j,k}^\epsilon}-\widehat  {c_{j,k}^\epsilon})^2
{\bf{1}}\{|\widetilde {c^\epsilon_{j,k}}| >
\lambda_n\}{\bf{1}}\{|\widehat  {c^\epsilon_{j,k}}| >
\frac{\lambda_n}{2}\}\right]
\\&& +
E \left[\sum_{j_n}^{J_n}\sum_{k,\epsilon} (\widehat
{c_{j,k}^\epsilon})^2 {\bf{1}}\{|\widetilde {c^\epsilon_{j,k}}| \leq
\lambda_n\}{\bf{1}}\{|\widehat  {c^\epsilon_{j,k}}| >
\frac{\lambda_n}{2}\}\right]\\&& + E
\left[\sum_{j_n}^{J_n}\sum_{k,\epsilon} (\widetilde
{c_{j,k}^\epsilon})^2 {\bf{1}}\{|\widetilde {c^\epsilon_{j,k}}| >
\lambda_n\}{\bf{1}}\{|\widehat {c^\epsilon_{j,k}}|
 \leq \frac{\lambda_n}{2}\}\right]
\\&=&T_1+T_2+T_3+T_4.
\end{eqnarray*}}
Using Proposition \ref{variance}, we have
\begin{eqnarray}\label{T1}
 T_1\leq
 K \: \frac{\log(n)}{n}
2^{dj_n} \leq K \: \frac{\left(\log(n)\right)^2}{n}.
\end{eqnarray}
 For the study of $T_2$,
we separate the cases where the wavelet coefficients are larger or
smaller than the thresholding level $\lambda_n/4$. By
Cauchy-Schwartz Inequality, we have

{\small \begin{eqnarray*}
T_2&=& E
\left[\sum_{j_n}^{J_n}\sum_{k,\epsilon} (\widetilde
{c_{j,k}^\epsilon}-\widehat  {c_{j,k}^\epsilon})^2
{\bf{1}}\{|\widetilde {c^\epsilon_{j,k}}| >
\lambda_n\}{\bf{1}}\{|\widehat  {c^\epsilon_{j,k}}| >
\frac{\lambda_n}{2}\}\left({\bf{1}}\{|c^\epsilon_{j,k}| \leq
\frac{\lambda_n}{4} \}+{\bf{1}}\{|c^\epsilon_{j,k}| >
\frac{\lambda_n}{4}\}\right)\right]
\\
&\leq & \sum_{j_n}^{J_n}\sum_{k,\epsilon} \left[E(\widetilde
{c_{j,k}^\epsilon}-\widehat  {c_{j,k}^\epsilon})^4\right]^{1/2}
\left[\mathbb{P}\left(|\widehat
{c_{j,k}^\epsilon}-c_{j,k}^\epsilon|
> \frac{\lambda_n}{4}\right)\right]^{1/2}\\&& +\sum_{j_n}^{J_n}\sum_{k,\epsilon}
E(\widetilde {c_{j,k}^\epsilon}-\widehat  {c_{j,k}^\epsilon})^2
{\bf{1}}\{|c^\epsilon_{j,k}| >  \frac{\lambda_n}{4}\}.
\end{eqnarray*}}
 Observe that, for any $j,k,\epsilon$, we have
\begin{eqnarray}\label{coefsup}
 |\widetilde
{c_{j,k}^\epsilon}|\vee |\widehat {c_{j,k}^\epsilon}|\leq
2^{j{d}/2}(\|\psi\|_\infty^{d}\vee \|\phi\|_\infty^{d}).
\end{eqnarray}
For any $\delta>0$, we use now the standard Bernstein Inequality to
obtain
\begin{eqnarray}\label{clas}
\mathbb{P}\left(|\widehat
{c_{j,k}^\epsilon}-c_{j,k}^\epsilon|
> \frac{\lambda_n}{4}\right)\leq K\,n^{-\delta}
\end{eqnarray}
which is valid for a choice of $\kappa$ large enough.
Let us now fix $r$ in $]0,2[$. Applying Proposition \ref{variance} and
using (\ref{coefsup}), it follows
 {\small
 \begin{eqnarray*}
 T_2
&\leq & K \ \sum_{j_n}^{J_n}\sum_{k,\epsilon}
2^{jd}\left[\mathbb{P}\left(|\widehat
{c_{j,k}^\epsilon}-c_{j,k}^\epsilon|
> \frac{\lambda_n}{4}\right)\right]^{1/2} +\sum_{j_n}^{J_n}\sum_{k,\epsilon}
E(\widetilde {c_{j,k}^\epsilon}-\widehat  {c_{j,k}^\epsilon})^2
{\bf{1}}\{|c^\epsilon_{j,k}| > \frac{\lambda_n}{4}\}\nonumber\\
&\leq& K \left ( 2^{2dJ_n }n^{-\delta/2}+ u_n \left[
\left(\frac{\lambda_n}{4}\right)^{r}
\sum_{j_n}^{J_n}\sum_{k,\epsilon} {\bf{1}}\{|c^\epsilon_{j,k}|
>  \frac{\lambda_n}{4}\}\right]\right )
\end{eqnarray*}}
for $$u_n=\left(\frac{\lambda_n}{4}\right)^{-r}
\,\frac{\log(n)}{n}.$$
 Similarly, we have
{\small \begin{eqnarray*}
T_3&\leq& E
\left[\sum_{j_n}^{J_n}\sum_{k,\epsilon} (\widehat
{c^\epsilon_{j,k}})^2{\bf{1}}\{|\widetilde {c^\epsilon_{j,k}}| \leq
\lambda_n\} {\bf{1}}\{|\widehat  {c^\epsilon_{j,k}}| >
\frac{\lambda_n}{2}\}\left({\bf{1}}\{|c^\epsilon_{j,k}| \leq
\frac{\lambda_n}{4} \}+{\bf{1}}\{|c^\epsilon_{j,k}| >
\frac{\lambda_n}{4} \}\right)\right]
\\
&\leq& E \left[\sum_{j_n}^{J_n}\sum_{k,\epsilon} (\widehat
{c^\epsilon_{j,k}})^2 {\bf{1}}\{|\widehat  {c^\epsilon_{j,k}}| >
\frac{\lambda_n}{2}\}{\bf{1}}\{|c^\epsilon_{j,k}| \leq
\frac{\lambda_n}{4} \}\right]\\&&+ E
\left[\sum_{j_n}^{J_n}\sum_{k,\epsilon} (\widehat
{c^\epsilon_{j,k}})^2{\bf{1}}\{|\widetilde {c^\epsilon_{j,k}}| \leq
\lambda_n\} {\bf{1}}\{|c^\epsilon_{j,k}| > \frac{\lambda_n}{4}
\}\right]
\\
&\leq&K\sum_{j_n}^{J_n}\sum_{k,\epsilon}
2^{dj}\mathbb{P}\left(|\widehat {c^\epsilon_{j,k}}-c^\epsilon_{j,k}|
> \frac{\lambda_n}{4}\right) +\left(\frac{\lambda_n}{4}\right)^r\sum_{j_n}^{J_n}\sum_{k,\epsilon}v_n
{\bf{1}}\{|c^\epsilon_{j,k}| >  \frac{\lambda_n}{4}\}
\end{eqnarray*}}
for
\begin{eqnarray*}
 v_n&=&2 \
\left(\frac{\lambda_n}{4}\right)^{-r}\left[E(\widetilde
{c_{j,k}}-\widehat {c_{j,k}})^2 +E(\widetilde
{c_{j,k}})^2{\bf{1}}\{|\widetilde {c_{j,k}}|
\leq  \lambda_n \}\right]\\
&\leq &2 \left( K u_n+ 4^{r} \lambda_n^{2-r}\right)
\end{eqnarray*}
implying that
{\small \begin{eqnarray*}\label{T3} T_3&\leq& K
\left(\frac{2^{2{dJ_n}}}{n^{\delta}}+
(u_n+\lambda_n^{2-r})\left[\left(\frac{\lambda_n}{4}\right)^r\sum_{j_n}^{J_n}\sum_{k,\epsilon}
{\bf{1}}\{|c^\epsilon_{j,k}|
>  \frac{\lambda_n}{4}\}\right]\right).
\end{eqnarray*}}
Using (\ref{coefsup}) and Proposition \ref{LDL}, we get
\begin{eqnarray*}\label{T4}
T_4 \leq  K
\sum_{j_n}^{J_n}\sum_{k,\epsilon}2^{dj}\mathbb{P}\left(|\widetilde
{c^\epsilon_{j,k}}-\widehat{c^\epsilon_{j,k}}|
> \frac{\lambda_n}{2}\right) \leq K \: 2^{2 dJ_n }\,n^{-\delta}.
\end{eqnarray*}
Combining the bounds of $T_1,T_2,T_3,T_4$  and choosing $j_n,J_n$ as
indicated in Theorem \ref{AdaptH}, we get for $\delta\geq 6$
\begin{eqnarray*} E\|\widetilde {c_{HL}}-c\|_2^2&\leq
& 2 \ E\|\widehat {c_{HL}}-c\|_2^2+K\,\rho_n
\end{eqnarray*}
where
 \begin{eqnarray*}
\rho_n &=& \frac{\left(\log(n)\right)^2}{n}
+ \left(\frac{\log n}{n}\right)^{1-\frac{r}{2}}
\left(\frac{\lambda_n}{4}\right)^r\sum_{j_n}^{J_n}\sum_{k,\epsilon}
{\bf{1}}\{|c^\epsilon_{j,k}|
> \frac{\lambda_n}{4}\}+ \frac{1}{n (\log(n))^2}.
\end{eqnarray*}
On the one hand, let us suppose that $c$ belongs to the weak Besov
$\mathcal{W}_{L}(\frac{2{d}}{2s+{d}})$ which means that (for
$r:={2d}/(2s+{d})$)
\begin{eqnarray*}
\left(\frac{\lambda_n}{4}\right)^r\sum_{j_n}^{J_n}\sum_{k,\epsilon}
{\bf{1}}\{|c^\epsilon_{j,k}|
>  \frac{\lambda_n}{4}\}& \leq &K.
\end{eqnarray*}
It follows that
\begin{eqnarray*}
\rho_n &\leq& K \:
\left(\frac{\log(n)}{n}\right)^{\frac{2s}{2s+{d}}}.
\end{eqnarray*}
Using the standard result given in Theorem \ref{maxisetstandard}
when direct observations are
 available, we also have
 \begin{eqnarray*}
E\|\widehat {c_{HL}}-c\|_2^2&\leq& K \:
\left(\frac{\log(n)}{n}\right)^{\frac{2s}{2s+{d}}}
\end{eqnarray*}
as soon as $c\in \mathcal{W}_{L}(\frac{2{d}}{2s+{d}})\cap
\mathcal{B}^s_{2\infty}$. This ends the proof of the first part
 of (\ref{maxih1}) of Theorem \ref{maxih}.

\subsubsection{Proof  of Equivalence $(\ref{maxih1})$ (second step: $\Longleftarrow$)}
Suppose that there exists $M$ such that for any $n$, $E\|\widetilde
{c_{HL}}-c\|_2^2\leq M \left(n^{-1}\,
\log(n)\right)^{\frac{2s}{2s+{d}}}.$ Since
$$\sum_{j>J_n}\sum_{k,\epsilon}(c_{j,k}^\epsilon)^2 \leq E\|\widetilde {c_{HL}}-c\|_2^2,$$
and choosing $J_n$ as indicated in Theorem \ref{AdaptH}, we obtain
$$\sum_{j>J_n}\sum_{k,\epsilon}(c_{j,k}^\epsilon)^2 \leq M
\left(\frac{\log(n)}{n}\right)^{\frac{2s}{2s+{d}}} \leq M
\left(2^{{d}(1-J_n)}\right)^{\frac{2s}{2s+{d}}} \leq K \left(
2^{-2J_n}\right)^{\frac{ds}{2s+{d}}}.$$ Using Definition \ref{Besov}
of the strong Besov spaces, we deduce that $c$ belongs necessarily
to $\mathcal{B}^{\frac{d s}{2s+{d}}}_{2 \infty}$. Let us now study
the sum of the square of the details when the details coefficients
are small
{\small{\begin{eqnarray*}\label{H0}
\sum_{j\geq
0}\sum_{k,\epsilon}(c_{j,k}^\epsilon)^2 {\bf{1}}\{|c^\epsilon_{j,k}|
\leq \frac{\lambda_n}{2} \}&=&
\left[\sum_{j<j_n}+\sum_{j=j_n}^{J_n}+\sum_{j>J_n}\right]
\left[\sum_{k,\epsilon}(c_{j,k}^\epsilon)^2
{\bf{1}}\{|c^\epsilon_{j,k}|
\leq \lambda_n/2 \}\right] \nonumber \\
&\leq &H_1+H_2+H_3.
\end{eqnarray*}}}
 Since we have already proved that
$c \in \mathcal{B}^{\frac{d s}{2s+{d}}}_{2 \infty}$ and taking
$\lambda_n$ as indicated in Theorem \ref{AdaptH}, we deduce
\begin{eqnarray*}\label{H3}
 H_3 &\leq& \sum_{j >
J_n}\sum_{k,\epsilon}(c_{j,k}^\epsilon)^2 \leq K   2^{-2J_n\frac{ d
s}{2s+{d}}} \leq  K \
\left(\frac{\lambda_n}{2}\right)^{\frac{4s}{2s+{d}}}.
\end{eqnarray*}
\noindent Taking $j_n$ as in Theorem \ref{AdaptH}, we get
\begin{eqnarray*}\label{H1}
H_1 &\leq&
K\sum_{j<j_n}2^{{d}j}\left(\frac{\lambda_n}{2}\right)^2 \leq K
\log(n) \left(\frac{\lambda_n}{2}\right)^2 \leq
K\;\left(\frac{\lambda_n}{2}\right)^{\frac{4s}{2s+{d}}}.
\end{eqnarray*}
\noindent  Observe that
\begin{eqnarray*}
H_2&=&E
\left[\sum_{j_n}^{J_n}\sum_{k,\epsilon}(c_{j,k}^\epsilon)^2
{\bf{1}}\{|c^\epsilon_{j,k}| \leq \frac{\lambda_n}{2} \}
\left({\bf{1}}\{|\widetilde {c^\epsilon_{j,k}}| \leq \lambda_n \}+
{\bf{1}}\{|\widetilde {c^\epsilon_{j,k}}| > \lambda_n
\}\right)\right].
\end{eqnarray*}
\noindent Remembering that
\begin{eqnarray*} E\left[
\sum_{j_n}^{J_n}\sum_{k,\epsilon}(c_{j,k}^\epsilon)^2{\bf{1}}\{|\widetilde
{c^\epsilon_{j,k}}| \leq \lambda_n \}\right] & \leq& E\|\widetilde
{c_{HL}}-c\|_2^2
\end{eqnarray*}
 and using Proposition \ref{LDL} and (\ref{clas}), we get
\begin{eqnarray*}\label{H2}
 H_2 &\leq&E\|\widetilde
{c_{_{HL}}}-c\|_2^2+
\sum_{j_n}^{J_n}\sum_{k,\epsilon}(c_{j,k}^\epsilon)^2
\mathbb{P}(|\widetilde {c^\epsilon_{j,k}}-\widehat{c^\epsilon_{j,k}}| >  \frac{\lambda_n}{4} )\nonumber \\
&& + \sum_{j_n}^{J_n}\sum_{k,\epsilon}(c_{j,k}^\epsilon)^2
\mathbb{P}(|\widehat {c^\epsilon_{j,k}}-c^\epsilon_{j,k}| >
 \frac{\lambda_n}{4} )\nonumber\\
 &\leq & M \ \left( \frac{\log(n)}{n}\right)^{\frac{2s}{2s+{d}}}+
K \ \|c\|_2^2\;n^{-\delta}  \leq  K
\left( \frac{\lambda_n}{2} \right)^{\frac{4s}{2s+{d}}}
\end{eqnarray*}
as soon as $\delta$ is larger than $1$. Combining  using Definition \ref{WLBesov} of the
local weak Besov space, we conclude that $c \in \mathcal{W}_{L}(r)$ with
$r$ such that $2-r=4s/(2s+{d})$. Hence, we end the proof of the
indirect direction of (\ref{maxih1}).

\subsection{Proofs of Proposition \ref{linclus} and Corollary
\ref{bestmaxi}}\label{preuvelienbesov} The proof of the {\bf {large
inclusion}} given in Proposition \ref{linclus} follows immediately
from the definitions of the functional spaces. Denote
$c_{j,k}^\epsilon$ the sequence of wavelet coefficients of a
function $c$. Since we have
\begin{eqnarray*} &&\sup_{0<\lambda \leq
1}\lambda^{r-2}\sum_{j\geq 0}\sum_{k,\epsilon}
 (c_{j,k}^\epsilon)^2{\bf{1}}\{|c_{j,k}^\epsilon| \leq
\lambda \}\\
&=&\sup_{0 < \lambda \leq 1}\lambda^{r-2}\sum_{j \geq
0}\sum_{k,\epsilon}
 (c_{j,k}^\epsilon)^2{\bf{1}}\{|c_{j,k}^\epsilon| \leq
\lambda \}\left[ {\bf{1}}\{\sum_k (c^\epsilon_{j,k})^2 \leq 2^{{d}j}
\lambda^2 \}\right.\\&&\left. \hspace{3cm}+ {\bf{1}}\{\sum_k
(c^\epsilon_{j,k})^2
>
2^{{d}j} \lambda^2 \}\right]\\
 &\leq& \sup_{0<\lambda \leq 1}\lambda^{r-2}\sum_{j\geq 0}\sum_{k,\epsilon}
 (c_{j,k}^\epsilon)^2{\bf{1}}\{\sum_k (c^\epsilon_{j,k})^2 \leq
2^{{d}j} \lambda^2 \}\\&& \hspace{3cm}+ K \: \sup_{0<\lambda \leq 1}\lambda^{r}\sum_{j\geq 0}2^{{d}j}\sum_{\epsilon}{\bf{1}}\{\sum_k
(c^\epsilon_{j,k})^2> 2^{{d}j} \lambda^2 \},
\end{eqnarray*}
it follows from
Definition \ref{WGBesov} that
$$c \in \mathcal{W}_{G}(r) \Rightarrow c \in
\mathcal{W}_{L}(r).$$
 To establish the {\bf {strict inclusions}}, we build a
sparse function belonging to $ \mathcal{B}^{\frac{ds}{2s+d}}_{2
\infty} \cap \mathcal{W}_{L}(\frac{2d}{2s+d})$ but not to
$\mathcal{W}_{G}(\frac{2{d}}{2s+d})$. Let us choose a real number
$\alpha$ such that $\frac{d}{2} \leq \alpha < s+\frac{d}{2}$. Let us
consider a function $c$ with the sparse sequence $c_{j,k}^\epsilon$
such that at each level $j \in \mathbb{N}$ and at each $\epsilon \in
S_d$, only $\lfloor2^{\frac{2d\alpha}{2s+d} j}\rfloor$ wavelet
coefficients take the value $(2^d-1)^{-1}2^{-alpha j}$. The others
are equal to $0$. For all $0<\lambda\leq 1,$ let $j_\lambda$ be such
that $2^{j_\lambda}=
\left((2^d-1)\lambda\right)^{-\frac{1}{\alpha}}$. We get
\begin{eqnarray*}
\sum_{j\geq
0}\sum_{k,\epsilon}{\bf{1}}\{|c_{j,k}^\epsilon|>\lambda\}&=&
\sum_{j<j_\lambda}\sum_{k,\epsilon}{\bf{1}}\{|c_{j,k}^\epsilon|>\lambda\}\\
&\leq& K \: 2^{\frac{2d\alpha}{2s+d} j_\lambda}\leq K \:
\lambda^{-\frac{2d}{2s+d}}.
\end{eqnarray*}
implying that $$\sup_{0<\lambda \leq
1}\lambda^{\frac{2{d}}{2s+d}}\sum_{j\geq
0}\sum_{k,\epsilon}{\bf{1}}\{|c_{j,k}^\epsilon|>\lambda\}<\infty,$$
 and the function $c$ belongs to the local weak Besov space
 $\mathcal{W}_{L}(\frac{2d}{2s+ d}).$
Next, put
 $\alpha^\prime=(4\alpha s + 2 sd + d^2)/(2(2s+d))$
 and observe that $\alpha^\prime<s+d/2$ since $\alpha<s+d/2$.
For all $0<\lambda \leq 1$ let $j_\lambda^*$ be such that $
2^{j_\lambda^{*}}=\left((2^d-1)\lambda\right)^{-\frac{1}{\alpha^\prime}}$.
We get
\begin{eqnarray*}
\sum_{j\geq
0}2^{dj}\sum_{\epsilon}{\bf{1}}\{\sum_k (c_{j,k}^\epsilon)^2>
2^{{d}j}\lambda^2\}&\geq& (2^{d}-1) \: \sum_{j<j_\lambda^*}2^{dj}\\
&\geq& 2^{{d} j_\lambda^*-1} \geq K \
\lambda^{-\frac{d}{\alpha^\prime}}.
\end{eqnarray*}
implying  $$\sup_{0<\lambda \leq
1}\lambda^{\frac{2{d}}{2s+{d}}}\sum_{j\geq
0}2^{{d}j}\sum_{\epsilon}{\bf{1}}\{\sum_k
(c_{j,k}^\epsilon)^2>2^{{d}j}\lambda^2\}= \infty.$$ It follows that
the function $c$ does not belong to the global weak Besov
$\mathcal{W}_{G}(\frac{2{d}}{2s+{d}})$ which ends the proof of
Proposition \ref{linclus}.  Notice that the function $c$ belongs to
the strong Besov space $\mathcal{B}^{\frac{ ds}{2s+d}}_{2 \infty} $
because for all $(j,\epsilon)$
\begin{eqnarray*}
\sum_{k,\epsilon}(c_{j,k}^\epsilon)^2& \leq &
2^{\frac{2d\alpha}{2s+d} j} 2^{-2 \alpha j} \leq
2^{-\frac{2ds}{2s+d}j}
\end{eqnarray*}
implying that $$\sup_{J \geq 0}2^{\frac{2ds}{2s+d}J}\sum_{j\geq
J}\sum_{k,\epsilon}(c_{j,k}^\epsilon)^2<\infty.$$
 Corollary \ref{bestmaxi} is also proved too.

\section{Appendix}
  This section aims at proving  (\ref{a}) of Proposition \ref{LDL}. In the sequel we fix the indices $j$
and $ k=(k_1,\ldots,k_{d})$ and take without loss of generality
$\epsilon=2^{d}-1$.
 For any  $i=1,\ldots, n$ (the observation index)  and any  $m=1,\ldots d$
 (the coordinate index), let us introduce the following notations
\begin{eqnarray*}
\Delta(X^m_i)&=&\widehat{F_m}(X^m_i)-F_m(X^m_i),\\
\xi_{j}(X^m_i)&=&\psi_{j,k_m}(\widehat{F_m}(X^m_i))-\psi_{j,k_m}(F_m(X^m_i)),\\
N_{j}(m)&=& \#\left\{i \in \{1, \dots, n\}; \xi_{j}(X^m_i)\not=
0\right\},\end{eqnarray*}
as univariate quantities and
\begin{eqnarray*}
 \xi_{j}(X^1_i,\ldots,X^m_i)&=&
\psi_{j,k}^\epsilon(\widehat{F_1}(X^1_i),\ldots,\widehat{F_{d}}(X^{d}_i)
)-\psi_{j,k}^\epsilon({F_1}(X^1_i),\ldots,{F_{d}}(X^{d}_i))\\
N_{j}&=& \#\left\{i \in \{1, \dots, n\};
\xi_{j}(X^1_i,\ldots,X^{d}_i)\not= 0\right\}
\end{eqnarray*}
 as $d-$variate quantities.
As previously remarked in Genest et al.  \cite{MasTrib} (for $d=2$),
from the definitions above we have
{\small{
\begin{eqnarray}\label{biuni}
\xi_{j}(X_i^1,\ldots,X_i^{d})&=&
\prod_{m=1}^{d}\xi_{j}(X_i^m)+\sum_{m_1=1}^{d}\left[
\psi_{j,k_{m_1}}^\epsilon({F_{m_1}}(X^{m_1}_i))
\prod_{\substack{m=1\\m\not=m_1}}^{d}\xi_{j}(X_i^{m}) \right]\nonumber \\
&&+\sum_{\substack{m_1,m_2=1\\m_1\not=m_2}} ^{d}\left[
\psi_{j,k_{m_1}}^\epsilon({F_{m_1}}(X^{m_1}_i))\psi_{j,k_{m_2}}^\epsilon({F_{m_2}}(X^{m_2}_i))
\prod_{\substack{m=1\\m\not=m_1,m_2}}^{d}\xi_{j}(X_i^{m}) \right]\nonumber\\
&&+\ldots +\sum_{m_1=1}^{d}\left[\xi_{j}(X_i^{m_1})
\prod_{\substack{m=1\\m\not=m_1}}^{d}
\psi_{j,k_{m}}^\epsilon({F_{m}}(X^{m}_i))\right].
\end{eqnarray}}}
 In the sequel, for
$m=1,\ldots ,d$, we denote by $T_{m,j}(X_i)$ any term of the type
$$
\left[\psi_{j,k_{1}}^\epsilon({F_{1}}(X^{1}_i))\times
\ldots\times\psi_{j,k_{d-m}}^\epsilon({F_{d-m}}(X^{d-m}_i))
\right]\;\left[\xi_{j}(X_i^{d-m+1})\times\ldots\times
\xi_{j}(X_i^{d})\right]
$$
i.e. such that there are exactly $m$ factors $\xi_j(X_i^{\cdot})$
appearing in the product. The cardinality of such terms $T_{m,j}(X_i)$
is equal to $C_{d}^m=\frac{d!}{m!(d-m)!}$. Observe that the number
of terms in (\ref{biuni}) is $2^{d}-1$. It is fundamental to notice
that there is no term $T_{0,j}(X_i)=\prod_{m=1,\ldots,d
}\psi_{j,k_{m}}^\epsilon({F_{m}}(X^{m}_i))$.
\subsection{Technical lemmas}
 We begin by giving technical lemmas.
\vspace{0.3cm}
\begin{lemma}\label{kief} There exists a universal
constant $K_0$ such that for any $m \in \{1, \dots,d\}$
\begin{eqnarray*}
 \forall t>0,\;
\mathbb{P}(\max_{1\leq i \leq n}|\Delta(X^m_i)|> t)&\leq& K_0
\exp(-2n t^2).
\end{eqnarray*}
\end{lemma}
Lemma \ref{kief} is a consequence of the Dvoreski-Kiefer-Wolfovitz inequality.
For the interested reader, the detailed proof of this lemma is given in Autin et al. \cite{ALT}.

\vspace{0.3cm}

\begin{lemma}\label{triboulnix} Let $\delta>0$ and $n$
be an integer such that $n \log(n) \geq 2 (\delta^{-1} \vee 1)$.
Then, there exists $K_1>0$ such that for any level $j$ satisfying
$$2^j\leq
\frac{1}{3}\,\left(\frac{2n}{\delta\log(n)}\right)^{1/2},$$
 and for any $m \in \{1, \dots,d\}$,
\vspace{0.3cm}\begin{eqnarray*}\label{nx}
\mathbb{P}(N_{j}(m)>(L+3)n2^{-j})\quad \vee \quad
\mathbb{P}(N_{j}>d(L+3)n2^{-j})&\leq & K_1 \;n^{-\delta}.
\end{eqnarray*}
\end{lemma}
 The proof of Lemma \ref{triboulnix}  can be found in  Autin et al. \cite{ALT}
\vspace{0.3cm}\begin{lemma}\label{LDZ} Let us assume that $c$
belongs to $L_\infty$ and let $(j,N)\in \N^2$. For all $1\leq p\leq
q\leq d$, for all subsets $\mathcal{S}_p$ and $\mathcal{S}_{q-p}$ of
$\{1,\dots,d\}$ with cardinalities equal to $p$ and $q-p$ having no
common component, let us put for $i=1\ldots,n$,
\vspace{0.3cm}\begin{eqnarray}\label{Z} Z_i(\mathcal{S}_p,
\mathcal{S}_{q-p})=\prod_{m \in \mathcal{S}_p}
\psi_{j,k_{m}}(F_{m}(X^{m}_i))\prod_{m^\prime \in \mathcal{S}_{q-p}}
(\psi^{(1)})_{j,k_{m^\prime}}(F_{m^\prime}(X^{m^\prime}_i)),
\end{eqnarray}
\noindent with the following notation $\psi_{j,k}^{(1)}(.)=2^{j/2}\psi'(2^j.-k)$.\\
\noindent For any $\mu\geq 2K_3 2^{-jq/2} $,
 we have
\vspace{0.3cm}\begin{eqnarray*}
\mathbb{P}\left(\left|\frac{1}{N}\sum_{i=1}^{N} Z_i(\mathcal{S}_p,
\mathcal{S}_{q-p})\right|>\mu \right)&\leq
&2\exp\left(-K_2N\left(\mu^2\wedge \mu2^{1-jq/2} \right) \right)
\end{eqnarray*}
 where $K_2,K_3$ are
constants such that {\footnotesize $$K_3\geq (L+1)^q \|c\|_\infty
\|\psi\|_\infty^{p}\|\psi'\|_\infty^{q-p}, \quad K_2\leq \frac{1}{8}
\|\psi\|_\infty^{-p}\|\psi'\|_\infty^{p-q}\left(K_3^{-1} \vee
6\right).$$}
\end{lemma}

 Lemma \ref{LDZ} is a direct application of the
Bernstein Inequality  with \vspace{0.3cm}

$\left|EZ_i(\mathcal{S}_p, \mathcal{S}_{q-p})\right|$

\begin{eqnarray*}&=&\left|\,\int_{[0,1]^{d}}
\prod_{m \in \mathcal{S}_p}\psi_{j,k_{m}}(u_{m})\prod_{m' \in
\mathcal{S}_{q-p}} \psi_{j,k_{m^\prime}}^{(1)}(u_{m^\prime})\:
c(u_1,\ldots,u_{d}) du_1\times \ldots \times
du_{d}\;\right|\\
&\leq& K_3\, 2^{-jq/2}
\end{eqnarray*} and in the same way,
\vspace{0.3cm}\begin{eqnarray*}
Var(Z_i(\mathcal{S}_p, \mathcal{S}_{q-p}))&\leq& (L+1)^q \|c\|_\infty \|\psi\|_\infty^{2p}\|\psi'\|_\infty^{2(q-p)}\\
\end{eqnarray*}
and
\begin{eqnarray*}
 \left|Z_i(\mathcal{S}_p, \mathcal{S}_{q-p})\right|&\leq& \|\psi\|_\infty^{p}\|\psi'\|_\infty^{q-p}2^{jq/2}.
\end{eqnarray*}

\subsection{Proof of Proposition \ref{LDL}}
By Equality (\ref{biuni}), we have for any $\lambda>0$
\vspace{0.3cm}\begin{eqnarray*}
\mathbb{P}\left(|\widehat{c^\epsilon_{j,k}}-\widetilde{c^\epsilon_{j,k}}|
>\lambda\right) &\leq& \sum_{m=1}^{d}
C_{d}^mL_m
\end{eqnarray*}
for \vspace{0.3cm}\begin{eqnarray*} L_m&=&
 \mathbb{P}\left(\left | \frac{1}{n}\sum_{i=1}^n T_{m,j}(X_i) \right |
>
\frac{\lambda}{2^{d}-1}\right).
\end{eqnarray*}
Using a Taylor expansion, the following inequality holds for
$i=1,\ldots, n$ and $m^\prime=1,\ldots,d $
\vspace{0.3cm}\begin{eqnarray*}\label{ksi2} |\xi_{j}(X^{m'}_i)| \leq
2^{j}|\Delta(X_i^{m'})|
(\psi^{(1)})_{j,k_{m^\prime}}(F_{m^\prime}(X^{m^\prime}_i))+
2^{\frac{3j}{2}-1}|\Delta(X_i^{m^\prime})|^2 \|\psi'\|_\infty
\end{eqnarray*}
implying that, for an associated  $\mathcal{S}_{d-m}$
{\small{\vspace{0.3cm}\begin{eqnarray*} |T_{m,j}(X_i)|&\leq&
\|\psi'\|_\infty^m  \mathop{\sum^{m}_{m'=0}}_{\mathcal{S}_{m-m'}\cap
\mathcal{S}_{d-m'}= \emptyset }
 2^{j(m+m'/2)}\left(\max_{m'=1,\ldots,m}
|\Delta(X_i^{m\prime})|\right)^{m+m^\prime}
\;\left|Z_i(\mathcal{S}_{d-m},\mathcal{S}_{m-m'})\right|.
\end{eqnarray*}}}
For $m=1,\ldots d $, let us introduce the following events {\small
\vspace{0.3cm}\begin{eqnarray*} {\cal D}_{0,m}=\left\{\max_{1 \leq i
\leq n} |\Delta(X^m_i)|\leq
\sqrt{\frac{\delta\log(n)}{2n}}\right\}&,&{\cal D}_{1,m}=\left\{
N_j(m) \leq
n_j=(L+3)n2^{-j}\right\},\\
{\cal D}_0=\bigcap_{m=1}^{d}{\cal D}_{0,m}&,&{\cal D}_1=\bigcap_{m=1}^{d}{\cal D}_{1,m}.
\end{eqnarray*}}
 It follows that for any $\mathcal{S}_p$ and any $\mathcal{S}_{q-p}$
\begin{eqnarray*} L_m&\leq&
 \mathbb{P}\left(\left(\left |
\frac{1}{n}\sum_{i=1}^n T_{m,j}(X_i) \right |
>
\frac{\lambda}{2^{d}-1}\right)\cap {\cal D}_0\cap {\cal
D}_1\right)+\mathbb{P}\left( {\cal D}_0^c\right)+\mathbb{P}\left(
{\cal D}_1^c\right)\\
&\leq&
 \sum_{m^\prime=0}^m\mathbb{P}\left( \left | \frac{1}{n_j}\sum_{i=1}^{n_j}
Z_i(\mathcal{S}_{d-m},\mathcal{S}_{m-m'}) \right |
>
\mu\right)+\mathbb{P}\left( {\cal D}_0^c\right)+\mathbb{P}\left(
{\cal D}_1^c\right)
\end{eqnarray*}
where $$\mu  = 2^{-j(m+m^\prime/2)}
\left(\frac{2n}{\delta\log(n)}\right)^{\frac{m+m^\prime}{2}}
\frac{2^j \|\psi'\|_\infty^{-m}(L+3)^{-1} \lambda}
{(2^{d}-1)(m+1)C^{\lfloor m/2 \rfloor}_m}.
$$
Fix $\kappa>0$ and take $\lambda= \sqrt{\frac{\kappa\log(n)}{n}}$.
 Using Lemma \ref{kief} and Lemma \ref{triboulnix}, one gets
\begin{eqnarray}\label{deltaxy} \mathbb{P}({\cal
D}_{0}^c)\vee \mathbb{P}({\cal D}_{1}^c) &\leq & d(K_0 \vee K_1)
\,n^{-\delta}
\end{eqnarray}
as soon as $2^j\leq
\frac{1}{3}\,\left(\frac{2n}{\delta\log(n)}\right)^{1/2}$. Since
$\mu\geq 2K_3\,2^{-j(d-m^\prime)/2}$, we apply Lemma \ref{LDZ} and
we obtain \begin{eqnarray*}
L_m&\leq&2\sum_{m^\prime=0}^m\exp\left[-K_22^{-j}n\left(\mu^2\wedge\mu2^{1-j(d
-m^\prime)/2}\right)\right] + d(K_0 \vee K_1)\,n^{-\delta}\leq
K\,n^{-\delta}
\end{eqnarray*}
as soon as \begin{eqnarray}\label{mu2} \mu\geq \left(
\frac{\delta}{K_2}\,\frac{2^j\log(n)}{n} \right)^{1/2} \vee \left(
\frac{\delta}{K_2}\,\frac{2^{j(2+d-m^\prime)/2}\log(n)}{2n} \right).
\end{eqnarray}
Let us restrict ourselves to the case where:
$$2^j\leq \left(\frac{n}{\log
n}\right)^{1/d}.$$ Assuming that $n$ is large enough and that
$\kappa$ is chosen large enough, Condition  (\ref{mu2}) on $\mu$ is
satisfied if, for any $m^\prime=0,\ldots m$
\begin{eqnarray*} d&\geq& \frac{2m+m^\prime-1}{m+m'}
\vee \frac{2m+d}{m+m'+1}.
\end{eqnarray*}
This condition  is always satisfied since $d \geq 2$. We obtain the announced result.

\noindent{\underline{Correspondances:}}

\noindent AUTIN Florent, LATP, universit\'e Aix-Marseille 1, Centre
de Math\'ematiques et Informatique, 39 rue F. Joliot Curie, 13453
Marseille Cedex 13 (autin@cmi.univ-mrs.fr).

\noindent LEPENNEC Erwan, LPMA, universit\'e Paris 7, 175 rue du
Chevaleret, 75013 Paris (lepennec@math.jussieu.fr).

\noindent TRIBOULEY Karine, LPMA, 175 rue du Chevaleret, 75013 Paris
and MODALX,  Universit\'e Paris 10-Nanterre, 200 avenue de la
R\'epublique, 92001 Nanterre Cedex (ktriboul@u-paris10.fr ).

\newpage

\newpage

\begin{table}
  \begin{center}
   \begin{tabular}[h]{|c|c||c||c|c|c|}
   \hline
  \multicolumn{2}{|c||}{Copula} &
  Method & \multicolumn{3}{|c|}{Boundaries Handling}
\\   \hline
   $c(\cdot)$&par. & & sym & per&ZeroPad \\\hline
  FGM & 1.0 &
Local&    0.007 (0.003)&0.079 (0.005)& 0.129 (0.010)\\
&&Block&    0.006 (0.002)&0.077 (0.008)& 0.141 (0.006)\\\hline
  normal& 0.0 &
Local&   0.002 (0.002)& 0.0004 (0.0004)&0.122 (0.005)\\
&&Block&    0.002 (0.002)& 0.0004 (0.0006)&0.105 (0.001)\\
  normal&$0.5$&
Local&   0.031 (0.007)&0.161 (0.011)& 0.179 (0.010)\\
&&Block&    0.032 (0.008)& 0.154 (0.011)& 0.202 (0.005)\\
 {\bf normal}&{\bf 0.9}&
Local&  {\bf  0.156 (0.011)}&{\bf 0.391 (0.008)}&{\bf 0.418 (0.006)}\\
&&Block& {\bf   0.140 (0.009)}&{\bf 0.381 (0.005)}&{\bf 0.491
(0.022)}\\\hline {\bf Student}&{\bf (0.5,1)}&
Local&  {\bf  0.326 (0.018)}&{\bf 0.460 (0.008)}&{\bf 0.544 (0.009)}\\
&&Block&  {\bf  0.324 (0.026)}&{\bf 0.458 (0.010)}&{\bf 0.585
(0.004)}\\\hline
  Clayton&$0.8$&
Local&   0.075 (0.013)&0.225 (0.010)& 0.252 (0.011)\\
&&Block&    0.095 (0.012)& 0.216 (0.011)&0.279 (0.005)\\\hline
  Frank&$4$&
Local&    0.021 (0.006)&0.149 (0.015)&0.212 (0.015)\\
&&Block&    0.013 (0.006)&0.134 (0.009)& 0.193 (0.006)\\\hline
 {\bf Gumbel}&{\bf 8.3}&
Local&  {\bf  0.701 (0.002)}&{\bf 0.849 (0.001)}&{\bf 0.866 (0.001)}\\
&&Block& {\bf   0.698 (0.002)}&{\bf 0.852 (0.001)}&{\bf 0.878 (0.001)}\\
  Gumbel&$1.25$&
Local&    0.038 (0.010)& 0.104 (0.005)& 0.172 (0.009)\\
&&Block&    0.052 (0.007)& 0.109 (0.004)& 0.173 (0.004)\\
  \hline
  \end{tabular}\end{center}
\caption{Relative $L_2$ estimation error for
$n=500$}\label{RE_2_n500}
\end{table}

\newpage

\begin{table}
\begin{center}
   \begin{tabular}[h]{|c|c||c||c|c|c|}
 \hline
  \multicolumn{2}{|c||}{Copula} &
  Method & \multicolumn{3}{|c|}{Boundaries Handling}
\\   \hline
   $c(\cdot)$&par. &  & sym & per&ZeroPad \\\hline
FGM & 1.0 &
Local& 0.0036 (0.0012)&0.0659 (0.0044)&0.0897 (0.0037)\\
&&Block&  0.0037 (0.0015)&0.0599 (0.0029)& 0.1068 (0.0041)\\\hline
normal&0.0 &
Local&  0.0006 (0.0005)& 0.0001 (0.0001)&0.0828 (0.0019)\\
&&Block& 0.0006 (0.0007)& 0.0001 (0.0001)& 0.0916 (0.0020)\\
normal&$0.5$&
Local& 0.0176 (0.0032)& 0.1449 (0.0040)& 0.1421 (0.0046)\\
&&Block&  0.0177 (0.0029)& 0.1329 (0.0036)&0.1518 (0.0055)\\
{\bf normal}&{\bf 0.9}&
Local& {\bf 0.1376 (0.0052)}&{\bf 0.3893 (0.0031)}&{\bf 0.4024 (0.0033)}\\
&&Block&{\bf  0.1330 (0.0045)}&{\bf 0.3813 (0.0027)}&{\bf 0.4261
(0.0046)}\\\hline {\bf Student}&{\bf (0.5,1)}&
Local&{\bf 0.2966 (0.0056)}&{\bf 0.4519 (0.0037)}&{\bf 0.5197 (0.0036)}\\
&&Block&{\bf 0.2881 (0.0058)}&{\bf 0.4467 (0.0029)}&{\bf 0.5230
(0.0028)}\\\hline Clayton&$0.8$&
Local& 0.0603 (0.0053)& 0.2073 (0.0046)&0.2127 (0.0041)\\
&&Block& 0.0596 (0.0054)& 0.1968 (0.0030)&0.2247 (0.0071)\\\hline
Frank&$4$&
Local& 0.01208 (0.0032)&0.1244 (0.0047)&0.1186 (0.0043)\\
&&Block&  0.0075 (0.0017)& 0.1137 (0.0035)& 0.1218 (0.0048)\\\hline
{\bf Gumbel}&{\bf 8.3}&
Local& {\bf 0.6975 (0.0015)}& {\bf 0.8511 (0.0004)}&{\bf 0.8664 (0.0003)}\\
&&Block& {\bf 0.6971 (0.0012)}&{\bf 0.8520 (0.0004)}&{\bf 0.8642 (0.0003)}\\
Gumbel&$1.25$&
Local&  0.0240 (0.0041)& 0.1022 (0.0030)&0.1392 (0.0029)\\
&&Block&  0.0336 (0.0042)&0.0988 (0.0026)& 0.1503 (0.0038)\\
\hline
\end{tabular}\end{center}
\caption{Relative $L_2$ estimation error for
$n=2000$}\label{RE_2_n2000}
\end{table}
\newpage

\begin{table}
\begin{center}
 \begin{tabular}[h]{|c|c||c||c|c|c|}
 \hline
  \multicolumn{2}{|c||}{Copula} &
  Method & \multicolumn{3}{|c|}{Empirical Loss Function}\\\hline
 $c(\cdot)$&par. & &$L_1$ & $L_2$ & $L_\infty$\\\hline
FGM & 1.0 &
Local&  0.062 (0.014)& 0.007 (0.003)& 0.189 (0.051) \\
&&Block& 0.061 (0.011)& 0.006 (0.002)& 0.175 (0.047) \\\hline
normal& 0.0&
Local& 0.038 (0.017)&0.002 (0.002)& 0.145 (0.062)\\
&&Block&0.038 (0.018)&0.002 (0.002)&0.129 (0.058) \\
normal&$0.5$&
Local&  0.118 (0.012)& 0.031 (0.007)& 0.539 (0.066) \\
&&Block& 0.112 (0.016)& 0.032 (0.008)& 0.555 (0.051) \\
{\bf normal}&{\bf 0.9}&
Local& {\bf 0.287 (0.026)}&{\bf 0.156 (0.011)}&{\bf 0.648 (0.020)} \\
&&Block&{\bf 0.205 (0.021)}&{\bf 0.140 (0.009)}&{\bf 0.644 (0.018)}
\\\hline {\bf Student}&{\bf (0.5,1)}&
Local&{\bf  0.290 (0.022)}&{\bf 0.326 (0.018)}&{\bf 0.791 (0.026) }\\
&&Block&{\bf 0.259 (0.018)}&{\bf 0.324 (0.026)}&{\bf 0.797 (0.035)}
\\\hline Clayton&$0.8$&
Local& 0.119 (0.014)& 0.075 (0.013)&0.658 (0.051) \\
&&Block&0.125 (0.018)& 0.095 (0.012)& 0.740 (0.040)\\\hline
Frank&$4$&
Local&  0.129 (0.017)& 0.021 (0.006)&0.329 (0.075) \\
&&Block& 0.092 (0.020)& 0.013 (0.006)& 0.321 (0.069) \\\hline {\bf
Gumbel}&{\bf 8.3}&
Local&{\bf  0.682 (0.015)}&{\bf 0.701 (0.002)}&{\bf 0.914 (0.001)}\\
&&Block&{\bf 0.629 (0.012)}&{\bf 0.698 (0.002)}&{\bf 0.915 (0.001)} \\
Gumbel&$1.25$&
Local& 0.099 (0.011)&0.038 (0.010)& 0.625 (0.104) \\
&&Block& 0.105 (0.012)&0.052 (0.007)& 0.749 (0.044) \\
\hline
\end{tabular}\end{center}
\caption{Relative $L_1$, $L_2$ and $L_\infty$ estimation errors for
$n=500$}\label{RE_1_2_infty_n500}
\end{table}

\newpage

\begin{table}
\begin{center}
 \begin{tabular}[h]{|c|c||c||c|c|c|}
 \hline
  \multicolumn{2}{|c||}{Copula} &
  Method & \multicolumn{3}{|c|}{Empirical Loss Function}\\\hline
 $c(\cdot)$&par.  & &$L_1$ & $L_2$ & $L_\infty$\\\hline
FGM & 1.0 &
Local& 0.0448 (0.00821) & 0.0036 (0.0012) &  0.1414 (0.0382) \\
&&Block& 0.04887 (0.0096)& 0.0037 (0.0015)& 0.1463 (0.0527)\\\hline
normal& 0.0 &
Local& 0.0181 (0.0087)&0.00063 (0.0005)&0.0673 (0.0332) \\
&&Block& 0.0190 (0.0092)& 0.0006 (0.0007)& 0.0669 (0.0284)
\\ normal&$0.5$&
Local& 0.0830 (0.0078)& 0.0176 (0.0032)& 0.4374 (0.0465) \\
&&Block& 0.0923 (0.0104)& 0.0177 (0.0029)& 0.4089 (0.0673)\\
{\bf normal}&{\bf 0.9}&
Local&{\bf 0.2048 (0.0160)}&{\bf 0.1376 (0.00522)}&{\bf 0.6400 (0.0114)} \\
&&Block&{\bf 0.1622 (0.0113)}&{\bf 0.1330 (0.0045)}&{\bf 0.6389
(0.0106)}
\\\hline {\bf Student}&{\bf (0.5,1)}&
Local& {\bf 0.2159 (0.0107)}&{\bf 0.2966 (0.0056)}&{\bf 0.7712 (0.0110)} \\
&&Block&{\bf 0.1955 (0.0095)}&{\bf 0.2881 (0.0058)}&{\bf 0.7669
(0.0133)}
\\\hline Clayton&$0.8$&
Local&  0.0862 (0.0068)& 0.0603 (0.0053)& 0.625 (0.0239) \\
&&Block& 0.1096 (0.0096)&0.0596 (0.0054)&0.6091 (0.0308) \\\hline
Frank&$4$&
Local& 0.0983 (0.0131)& 0.01208 (0.0032)&0.2635 (0.0569) \\
&&Block&0.0702 (0.0096)& 0.0075 (0.0017)& 0.2508 (0.0608)
\\\hline {\bf Gumbel}&{\bf 8.3}&
Local&  {\bf 0.6283 (0.0086)}&{\bf 0.6975 (0.0015)}&{\bf 0.9145 (0.0009)}\\
&&Block&{\bf 0.6223 (0.0058)}&{\bf 0.6971 (0.0012)}&{\bf 0.9143 (0.0007) }\\
Gumbel&$1.25$&
Local&  0.0720 (0.0075)&0.0240 (0.0041)& 0.5377 (0.0568) \\
&&Block& 0.0721 (0.0085)& 0.0336 (0.0042)& 0.6688 (0.0421) \\
\hline
\end{tabular}\end{center}
\caption{Relative $L_1$, $L_2$ and $L_\infty$ estimation errors for
$n=2000$}\label{RE_1_2_infty_n2000}
\end{table}

\newpage

\begin{table}[h]
\begin{center}
\begin{tabular}{|c|c||c|c||c|c||c|c|}\hline
&&$\hat \theta_1$&$E_1$&$\hat
\theta_2$&$E_2$&$\hat\theta_\infty$&$E_\infty$
\\\hline
Gaussian&Block&
 -0.01 &  0.0068  &  -0.01 &0.0001 &  -0.01
  &     0.0449\\
       Gaussian&Local&  -0.01&    0.0080 &  -0.01&   0.0002 &   0.01& 0.0847 \\\hline
Student&Block&
 (-0.11,91)&  0.0640 & (-0.11,91)&  0.0103 &   (-0.11,91)
    &    0.6639
\\
Student&Local &
 (0.07,40)&    0.0226 &   (0.07,40)&    0.0010  &  (0.02,100)&    0.1279
\\\hline
 Clayton &Block&  0.01&  0.0125&  0.01 &      0.0002&
0.01 &
      0.0395
\\
    Clayton &Local&0.01&    0.0135  &  0.01&   0.0004 &   0.01&  0.0942 \\\hline
    Frank&Block& 0.01&   0.0103& 0.01& 0.0002  &  0.01
   &      0.0467
\\
   Frank &Local&  0.01&    0.0115 &   0.01&  0.0003   & 0.07&   0.0825\\\hline
    Gumbel &Block&   1.00&  0.0093  &  1.00&   0.0002 &  1.00
   &   0.0462\\
           Gumbel &Local&  1.00&  0.0106 &   1.00&   0.0003 &   1.00&    0.0963\\\hline\hline
 All&Block&-0.01&Gaussian & -0.01&Gaussian &0.01&Clayton
 \\
&&&0.68\%& & 0.01 \%&   &   4.28   \%\\
All&Local&-0.01&Gaussian &-0.01&Gaussian &0.07&Frank
 \\
&&&  0.79  \%& &  0.02 \%&   &    7.98 \%\\
 \hline
\end{tabular}
\end{center}
\caption{{\small { Brent/Cac: distances between the benchmarks and
the parametrical families}}\label{Table_brent/cac}}
\end{table}

\newpage

\begin{table}[h]
\begin{center}
\begin{tabular}{|c|c||c|c||c|c||c|c|}\hline
&&$\hat \theta_1$&$E_1$&$\hat
\theta_2$&$E_2$&$\hat\theta_\infty$&$E_\infty$
\\\hline\hline
Gaussian &Block& 0.15&0.0396 &
0.14  &  0.0030&    0.10& 0.1337\\
Gaussian&Local&  0.14&    0.0492  &  0.13&   0.0041  &  0.10&
0.1437\\\hline
 Student& Block&(0.14,37)& 0.0376 &    (0.13,81)  & 0.0030  &
(0.08,61)
  &   0.1329\\
Student&Local & (0.14,95)&   0.0491  & (0.13,95)&    0.0041 & (0.09,
80)& 0.1411
\\\hline
 Clayton &Block& 0.15 & 0.0706 &   0.12  &0.0099 &  0.05
      &   0.1879\\
       Clayton &Local&  0.14&   0.0799 &   0.11&   0.0109   & 0.05&   0.1967\\\hline
Frank& Block&0.76  &  0.0301 &  0.83 & 0.0017 &   0.85
    &   0.0957\\
     Frank &Local& 0.75&    0.0393  &  0.80&   0.0027 &   0.54&  0.1355\\\hline
    Gumbel &Block& 1.10 & 0.0436  &   1.07  &  0.0069  & 1.02
    &  0.2309\\
     Gumbel &Local& 1.10&  0.0529 &   1.06&  0.0076  &  1.02&  0.2298 \\\hline\hline
    All&Block&0.76  &Frank &  0.83 &Frank&  0.85&Frank\\
&&&3.01   \%& &0.17    \%&   &
 6.61\%\\\hline
 All&Local& 0.75 &Frank &0.80 &Frank&0.54&Frank\\
&& & 3.93    \%& &   0.27   \%&   &
 10.64\%\\
  \hline
\end{tabular}
\end{center}
\caption{{\small {Brent/Exonmobil: distances between the benchmarks
and the parametrical families}}\label{Table_brent_exonmobil}}
\end{table}

\newpage

\begin{table}[h]\begin{center}
\begin{tabular}{|c|c||c|c||c|c||c|c|}\hline
&&$\hat \theta_1$&$E_1$&$\hat
\theta_2$&$E_2$&$\hat\theta_\infty$&$E_\infty$
\\\hline
Gaussian&Block& -0.11 &0.0233 &  -0.10&  0.0010 &  -0.07&
  0.0765
 \\
      Gaussian&Local&-0.11&    0.0243  & -0.10&    0.0011 &  -0.07& 0.0765 \\\hline
      Student&Block&(-0.11,61)&  0.0233&(-0.10,61)&  0.0011 &  (-0.06,61)
   &      0.0859
\\Student&Local & (-0.11,80)    &0.0239 &  (-0.10,80)  &  0.0011 &  (-0.06,63)&    0.0859\\\hline
Clayton &Block&0.01&  0.0801 &  0.01& 0.0104   &   0.01
   &   0.2924
\\
    Clayton &Local&   0.01 &   0.0805  &  0.01&   0.0105 &   0.01&    0.2924\\\hline
    Frank&Block& -0.57&  0.0148 &  -0.56&  0.0003 &  -0.50&     0.0456
\\
   Frank &Local& -0.58&    0.0155 &  -0.57&   0.0004 &  -0.48&   0.0433\\\hline
      Gumbel&Block&1.00&  0.0755& 1.00&  0.0090 &   1.00&
        0.2316
\\
Gumbel &Local&   1.00&   0.0760  &  1.00&  0.0092  &  1.00& 0.2316
\\\hline\hline All&Block&-0.57&Frank& -0.56 &Frank& -0.50
&Frank
 \\
&&&   1.48  \%&   &     0.03   \%&   &    3.69\%\\

All&Local&-0.58&Frank& -0.57&Frank& -0.48&Frank
 \\ &&&  1.54      \%&&  0.03      \%&   &     3.53\%\\
 \hline
\end{tabular}
\end{center}
\caption{{\small { DowJones/Oncedor: distances between the
benchmarks and the parametrical
families}}\label{Table_dowjones/oncedor}}
\end{table}

\newpage

\begin{table}[h]
\begin{center}
\begin{tabular}{|c|c||c|c||c|c||c|c|}\hline
&&$\hat \theta_1$&$E_1$&$\hat
\theta_2$&$E_2$&$\hat\theta_\infty$&$E_\infty$
\\\hline
Gaussian&Block&  0.30 & 0.0976 &   0.33& 0.0202 &   0.20&    0.4191
 \\
      Gaussian&Local& 0.26&   0.0699 &   0.32&    0.0234  &  0.11&    0.2785 \\\hline
      Student& Block& (0.28,8)&  0.0755  &  (0.29,8)&   0.0127& (0.18,11)&    0.3027
\\
Student&Local & (0.17,12)&    0.0846   & (0.17,6)& 0.0265&
(0.12,20)&    0.3748
\\\hline Clayton &Block&  0.40&   0.1064 &  0.36&   0.0318
& 0.26  &  0.4565
\\
    Clayton &Local&0.31&    0.0978 &   0.33&    0.0401 &   0.11&  0.3465 \\\hline
    Frank& Block& 1.58& 0.1094 &  1.88&  0.0333  &   0.57&    0.4366
\\
   Frank &Local&  1.38&    0.0687  &  1.73&    0.0401   & 0.79&  0.2762\\\hline
       Gumbel&Block&    1.19&   0.1081 &  1.17&    0.0414 & 1.09 &    0.4427
\\
 Gumbel &Local& 1.18&   0.0782 &   1.18&  0.0282  &  1.06&   0.3866 \\\hline\hline
All&Block&(0.28,8)&Student&(0.29,8)&Student&(0.18,11)&Student\\
&&&7.55 \%&   &  1.15 \%&   &    10.62\%
 \\
 All&Local&1.38&Frank&0.32&Gaussian&0.79&Frank\\
&&&  6.86  \%&   &  2.12 \%&   &   19.56  \%
 \\\hline
\end{tabular}
\end{center}
\caption{{\small { DowJones/Ftse100uk: distances between the
benchmarks and the parametrical
families}}\label{Table_dowjones/ftse100uk}}
\end{table}

\end{document}